\documentclass[11pt,reqno]{amsart}
\usepackage{graphicx}
\usepackage{color}
\usepackage{amsmath,amssymb,amsthm,amsfonts}
\usepackage{hyperref}
\usepackage{mathdots}
\usepackage[foot]{amsaddr}
\usepackage{cases}
\usepackage[numbers,sort&compress]{natbib}
\usepackage{multicol}
\usepackage{enumerate}
\usepackage{enumitem}
\usepackage{tikz}
\usepackage{tikz-3dplot}
\usetikzlibrary{calc}
\usepackage{caption}
\usepackage{subcaption}

\hypersetup{CJKbookmarks=true}

\hypersetup{hidelinks}
\numberwithin{equation}{section}

\def\Div{\operatorname{div}}

\newtheorem{thm}{Theorem}[section]

\newtheorem{lemma}[thm]{Lemma}
\newtheorem{theorem}[thm]{Theorem}
\newtheorem{proposition}[thm]{Proposition}

\theoremstyle{definition}
\newtheorem{definition}{Definition}[section]
\newtheorem{remark}{Remark}[section]

\begin{document}
	\title[]{Global existence and decay estimates for classical solutions to 2D compressible micropolar fluid equations with vacuum}
	\author{Guo houbin$^{1}$, ~~ Zhu changjiang$^{2,*}$, ~~ Zhu qiaolong$^{3}$}
    \thanks{*Corresponding author.}
	\address{$^{1,2,3}$School of Mathematics, South China University of Technology, Guangzhou 510640, China}
	\email{$^1$maghb@mail.scut.edu.cn}
	\email{$^2$machjzhu@scut.edu.cn}
	\email{$^3$mazhuqiaolong@mail.scut.edu.cn}

\begin{abstract}
In this paper, we consider the Cauchy problem for the two-dimensional compressible micropolar fluid equations with far field vacuum. Under the assumption that the initial energy is sufficiently small, we establish the global existence and uniqueness of strong and classical solutions, allowing for large oscillations in the initial data. The key to the proof is to derive {\it a priori} decay estimates for the pressure, the velocity gradient, and the microrotational velocity, which enable us to close the bootstrap argument. We also handle the difficulties arising from the coupling between the velocity field and the microrotational velocity.
\vspace{6mm}


\vspace{3mm}
\end{abstract}
\maketitle

\section{Introduction}
The micropolar fluid model, which is used to describe microrotation and spin inertia, is an important generalization of the Navier-Stokes equations. This model has been widely used to study complex fluids such as blood, polymer suspensions, and liquid crystals, cf. \cite{Condiff,Eringen,Lukaszewicz}. In this paper, we consider the two-dimensional micropolar fluid equations, which can be written as
\begin{align}\label{model}
	\left\{
	\begin{array}{l}
		\rho_t+\Div(\rho u)=0,
		\\
		(\rho u)_t+\Div(\rho u\otimes u)+\nabla P(\rho)=(\mu+\kappa)\Delta u+(\mu+\lambda-\kappa)\nabla\Div u+2\kappa\nabla^\perp\omega,
		\\
    	(\rho\omega)_t+\Div(\rho u\omega)+4\kappa\omega=\nu\Delta\omega+2\kappa\text{curl}~u,
	\end{array}
	\right.
\end{align}
where $\rho=\rho(x_1,x_2,t)$,~$u=(u_1(x_1,x_2,t),u_2(x_1,x_2,t),0)$, $\omega=(0,0,\omega(x_1,x_2,t))$, $P(\rho)$ represent the density,
velocity, microrotational velocity, and pressure, respectively. The viscosity coefficients $\mu$, $\lambda$, $\kappa$, $\nu$ satisfy the following physical constraints:
\begin{align}\label{viscosity}
\mu>0, \quad \kappa>0, \quad \nu>0, \quad \mu+\lambda-\kappa\geq0.
\end{align}
Here $\nabla^\perp\omega=(\partial_2\omega,-\partial_1\omega)$, $\text{curl}~u=\partial_1u_2-\partial_2u_1$. The pressure is given by $P(\rho)=A\rho^\gamma$, where $\gamma>1$ is the adiabatic exponent, and $A>0$ is a constant. Without loss of generality, we choose $A=1$. We consider the Cauchy problem for \eqref{model} with the initial data
\begin{align}\label{initial}
(\rho, \rho u, \rho\omega)(x,0)=(\rho_0, \rho_0u_0, \rho_0\omega_0)(x), \quad x\in\mathbb{R}^2,
\end{align}
and assume that $(\rho, u, \omega)$ vanish at infinity in some weak sense.

When microrotational effects are neglected ($\omega=\text{constant},~\kappa=0$), the model reduces to the compressible Navier-Stokes equations, which have been extensively studied. We briefly review some relevant results below. For the multi-dimensional case, in the absence of vacuum, the local existence and uniqueness of classical solutions were established in \cite{Nash, Serrin}. In \cite{Matsumura_Nishida}, Matsumura and Nishida first proved the global existence and uniqueness of classical solutions under the assumption that the initial data are sufficiently close to a nonvacuum equilibrium state in a suitable Sobolev space $H^s$. In \cite{Hoff1, Hoff2}, Hoff studied the existence of global weak solutions for discontinuous initial data with small energy. In the presence of vacuum, the local existence and uniqueness of strong or classical solutions were established in \cite{Li_Liang, Cho_Choe_Kim, Cho_Kim, Choe_Kim}. In \cite{Lions1}, Lions first established the existence of global weak solutions with large initial data in $\mathbb{R}^N$ when $\gamma\geq\frac{3N}{N+2}$, $N=2,3$. Later, Feireisl and his collaborators \cite{Feireisl} extended this result to the cases $\gamma>\frac{3}{2}$ for $N=3$. Huang, Li, and Xin \cite{Huang_Li_Xin} investigated the Cauchy problem for the three-dimensional isentropic Navier-Stokes equations with vacuum and showed the global existence and uniqueness of classical solutions with small total energy but possibly large oscillations. Subsequently, Li and Xin \cite{Li_Xin} established the global existence and uniqueness of solutions with small initial energy in two dimensions and obtained their large time asymptotic behavior. When $\gamma$ is sufficiently close to $1$, Hong, Hou, Peng, and Zhu \cite{Hong} proved the global existence and uniqueness of classical solutions with large initial energy.

The micropolar fluid equations have also been extensively studied. For the one-dimensional case, we refer to \cite{Mujakovic1, Mujakovic2, Mujakovic3, Mujakovic4, Mujakovic5, Chen, Feng_Zhu} for some results relevant to the present work. In the three dimensional case without vacuum, Liu and Zhang \cite{Liu_Zhang, Liu_Zhang1} investigated small perturbations around equilibrium states for the isentropic compressible micropolar fluid system, treating the cases without and with external forces, respectively. They established the global existence of smooth solutions and obtained optimal decay rates in $L^2$. In the presence of vacuum, Chen, Xu, and Zhang \cite{Chen_Xu} obtained the global existence of weak solutions for discontinuous initial data under the assumption of sufficiently small initial energy.  In the spirit of the work of Hong et al. \cite{Hong}, Hou and Peng \cite{Hou_Peng} proved the global existence and uniqueness of classical solutions in $\mathbb{R}^3$ with far field vacuum when $(\gamma-1)^{\frac{1}{3}}E_0$ is  sufficiently small. Their result permits large initial energy as $\gamma$ approaches $1$.
Inspired by \cite{Huang_Li_Xin}, Zhu and Tao \cite{Zhu_Tao} showed the global existence and uniqueness of classical solutions to the three-dimensional isentropic micropolar fluid system with small initial energy. They allow any fixed nonnegative far field density $\tilde\rho$, and assume $(\rho_0-\tilde\rho, P(\rho_0)-P(\tilde\rho))\in H^2\cap W^{2,q}$, with $3<q<6$.

Compared with the three-dimensional case, fewer results are available on the well-posedness of the two-dimensional compressible micropolar fluid equations. Wan and Zhang \cite{Wan_Zhang} established the global existence of classical solutions on $\mathbb{T}^2$ for arbitrarily large initial data allowing vacuum, assuming $\lambda=\rho^\beta$ with viscosity exponent $\beta>\frac{3}{2}$ and $1<\gamma<4\beta-3$. This density dependence introduces a term of order $\rho^\beta$ to  dominate the error terms arising from the effective viscous flux estimates, which yields a uniform-in-time upper bound for the density. With constant viscosities on $\mathbb{T}^2$,  Du, Yang, and You \cite{Du} established global existence for $\gamma>1$ under a sufficiently large bulk viscosity $b=2\mu+\lambda$  and a bound on $b^{\frac{1}{2}}\|\Div u_0\|_{L^2}$ uniform in $b$. In their argument, large $b$ strengthens the dissipation of $\Div u$ and makes the effective viscous flux contribution in the density estimate sufficiently small. For the Cauchy problem on $\mathbb{R}^2$, Lin, Liu, and Zhang \cite{Lin_Liu_Zhang} established the global existence and uniqueness of strong solutions for large initial data with density bounded away from vacuum. Their analysis also assumes a positive far field density, a sufficiently large bulk viscosity, and a uniform-in-$b$ bound on $b^{\frac{1}{2}}\|\Div u_0\|_{L^2}$. They further justified the combined limits of large bulk viscosity and vanishing angular viscosity.

The present work focuses on the Cauchy problem on $\mathbb{R}^2$ with finite mass and far-field vacuum. The viscosity coefficients are constant and  satisfy \eqref{viscosity}. We obtain global strong solutions, and classical solutions under additional regularity and compatibility conditions, without a density
dependent bulk viscosity, a largeness assumption on the bulk viscosity, or a preparation condition on $\Div u_0$.


Before stating the main results, we introduce the notation and conventions used throughout this paper. For $R>0$, set
\begin{align*}
B_R\triangleq\{x\in\mathbb{R}^2| |x|<R\}, \quad \int_{\mathbb{R}^2}fdx\triangleq\int f dx, \quad \dot{f}\triangleq f_t+u\cdot\nabla{f}.
\end{align*}
For $1\leq p\leq\infty$, $k\geq1$ and $\beta>0$, the standard homogeneous and inhomogeneous Sobolev spaces are defined as follows:
\begin{align*}
	\left\{
\begin{array}{l}
L^p=L^p(\mathbb{R}^2),\quad D^{k,p}=D^{k,p}(\mathbb{R}^2)=\{v\in L_{loc}^1(\mathbb{R}^2)|\nabla^kv\in L^p(\mathbb{R}^2)\},
\\
D^1=D^{1,2}, \quad W^{k,p}=W^{k,p}(\mathbb{R}^2), \quad H^k=W^{k,2},
\\
\dot{H}^\beta=\left\{f:\mathbb{R}^2\to\mathbb{R}\middle| \|f\|^2_{\dot{H}^\beta}=\int|\xi|^{2\beta}|\hat f(\xi)|^2d\xi<\infty\right\},
\end{array}
	\right.
\end{align*}
where $\hat f$ is the Fourier transform of $f$. Next, we give the definition of a strong solution to system \eqref{model} as follows:
\begin{definition}
	If all derivatives involved in \eqref{model} for $(\rho, u, \omega)$ are regular distributions, and equations \eqref{model} hold almost everywhere in $\mathbb{R}^2\times(0, T)$, then  $(\rho, u, \omega)$  is called a strong solution to \eqref{model}.
\end{definition}
In this paper, the initial total energy is defined as
\begin{align}
E_0=\int\left(\frac{1}{2}\rho_0|u_0|^2+\frac{1}{2}\rho_0|\omega_0|^2+\frac{1}{\gamma-1}P(\rho_0)\right)dx.\notag
\end{align}
The weight function is defined by
\begin{align}\label{weight}
	\bar{x}\triangleq(e+|x|^2)^ \frac{1}{2}\log^2(e+|x|^2).
\end{align}
We also introduce the effective viscous flux and the modified vorticity:
\begin{align}\label{flux_curl}
F=(2\mu+\lambda)\Div u-P(\rho),\quad \bar\omega=\text{curl}~u-\frac{2\kappa}{\mu+\kappa}\omega.
\end{align}
Furthermore, without loss of generality, we assume that there exists a positive constant $N_0$ such that the initial density $\rho_0$ satisfies
\begin{align}\label{tm}
 \frac{1}{2}\leq\int_{B_{N_0}}\rho_0dx\leq\int_{\mathbb{R}^2}\rho_0 dx=1.
\end{align}

The main results of this paper can be stated as follows:
\begin{theorem}\label{thm1.1}
	Let $M>0$, $\bar\rho\geq1$, $a>1$, $q>2$ and $\beta\in(0,1]$ be given numbers. Suppose that the initial data $(\rho_0, u_0, \omega_0)$ satisfy, in addition to \eqref{tm},
\begin{align}\label{initial1}
 \bar{x}^a\rho_0\in L^1\cap H^1\cap W^{1,q},\quad u_0\in\dot{H}^\beta\cap D^1, \quad \omega_0\in H^1,\quad \rho_0^{\frac{1}{2}}u_0,~\rho_0^{\frac{1}{2}}\omega_0\in L^2,
\end{align}
and
\begin{align}\label{initial2}
0\leq\rho_0\leq\bar\rho,\quad \|u_0\|_{\dot{H}^\beta}+\|\bar{x}^a\rho_0\|_{L^1}+\|\omega_0\|_{H^1}\leq M.
\end{align}
Then there exists a positive constant $\varepsilon$ depending only on $\mu$, $\lambda$, $\nu$, $\kappa$, $\gamma$, $a$, $\bar\rho$, $\beta$, $N_0$ and $M$ such that if
\begin{align}
E_0\leq\varepsilon,\notag
\end{align}
the problem \eqref{model}$-$\eqref{initial} admits a unique global strong solution $(\rho, u, \omega)$ satisfying,  for any finite $T>0$,
the following estimates and regularity properties:
\begin{align}\label{bound}
0\leq\rho(x,t)\leq2\bar\rho, \quad (x,t)\in\mathbb{R}^2\times[0,T],
\end{align}
\begin{align}\label{regularity}
\left\{
	\begin{array}{l}
\rho\in C([0,T];L^1\cap H^1\cap W^{1,q}),
\\
\bar{x}^a\rho\in L^\infty(0,T; L^1\cap H^1\cap W^{1,q}),
\\
\sqrt{\rho}u, \sqrt{\rho}\omega, \nabla{u}, \nabla{\omega}, \bar{x}^{-1}u, \omega, \sqrt{t}\sqrt{\rho}u_t, \sqrt{t}\sqrt{\rho}\omega_t\in L^\infty(0,T; L^2),
\\
\nabla{u}, \nabla{\omega}\in L^2(0,T; H^1)\cap L^{\frac{q+1}{q}}(0,T;W^{1,q}),
\\
\sqrt{t}\nabla{u}, \sqrt{t}\nabla{\omega}\in L^2(0,T; W^{1,q}),
\\
\sqrt{\rho}u_t, \sqrt{\rho}\omega_t, \sqrt{t}\nabla{u}_t, \sqrt{t}\nabla{\omega}_t, \sqrt{t}\omega_t, \sqrt{t}\bar{x}^{-1}u_t, \sqrt{t}\bar{x}^{-1}\omega_t\in L^2(\mathbb{R}^2\times(0,T)),
\end{array}
\right.
\end{align}
and for any $\alpha>1$, there exists a positive constant $N_1$ depending only on $\alpha$, $N_0$ and $M$ such that
\begin{align}\label{low_bound}
\inf_{0\leq t\leq T}\int_{B_{N_1(1+t)\log^{\alpha}(e+t)}}\rho(x,t)dx\geq\frac{1}{4}.
\end{align}
Moreover, for $t\geq1$, $(\rho,u,\omega)$ has the following decay estimates:
\begin{align}\label{decay}
\left\{
\begin{array}{l}
\|\nabla{u}(\cdot,t)\|_{L^p}+\|\omega(\cdot,t)\|_{L^p}\leq C(p)t^{-1+\frac{1}{p}}, \quad p\in[2,\infty),
\\
\|\nabla{\omega}(\cdot,t)\|_{L^p}\leq C(p)t^{-1}, \quad p\in[2,\infty),
\\
\|P(\cdot,t)\|_{L^r}\leq C(r)t^{-1+\frac{1}{r}}, \quad r\in(1,\infty),
\\
\|\nabla{F}(\cdot,t)\|_{L^2}+\|\nabla\bar\omega\|_{L^2}\leq Ct^{-1}.
\end{array}
\right.
\end{align}
\end{theorem}

\begin{remark}
The  condition $\bar{x}^a\rho_0\in L^1\cap H^1\cap W^{1,q}$ in \eqref{initial1}  is a sufficient technical hypothesis for the present proof, not a constitutive restriction of the micropolar model. Li and Xin \cite{Li_Xin} imposed the same weighted assumption.  The weighted $L^1$ condition restricts the mass at large distances. The weighted $ H^1\cap W^{1,q}$ regularity also restricts spatial variations: a smooth compactly supported density satisfies the condition, while a smooth profile with finite weighted mass can violate it through increasingly rapid oscillations at infinity. The weighted gradient norms need only be finite. \eqref{weight_rho_high_estimate} propagates the $H^1$ and $W^{1,q}$ bounds on $\bar x^a\rho$ required in \eqref{regularity}. Finite weighted mass alone does not provide the corresponding gradient bounds. Retaining only $\bar x^a\rho_0\in L^1$, together with suitable unweighted density regularity, would require a weaker weighted solution class and an adapted local existence and continuation argument.
\end{remark}

\begin{remark}
As in Li and Xin \cite{Li_Xin}, the decay estimates in \eqref{decay} are used to control the contribution of the effective viscous flux to the density equation at large times. This allows us to apply Zlotnik's inequality to obtain a uniform-in-time upper bound for the density and close the bootstrap argument. The optimality of these decay rates is not addressed here.
\end{remark}
If the initial data $(\rho_0, u_0, \omega_0)$ satisfy some additional  regularity and compatibility conditions, the global strong solution derived by Theorem \ref{thm1.1} becomes classical, as stated in the following theorem.
\begin{theorem}\label{thm1.2}
 In addition to the assumptions of Theorem \ref{thm1.1}, suppose further that the initial data $(\rho_0, u_0, \omega_0)$ satisfy
\begin{align}\label{ad_initial}
\nabla^2\rho_0,~\nabla^2P(\rho_0)\in L^2\cap L^q,\quad \bar{x}^{\delta_0}\nabla^2\rho_0,~\bar{x}^{\delta_0}\nabla^2P(\rho_0),~\nabla^2u_0,~\nabla^2\omega_0\in L^2,
\end{align}
for some constant $\delta_0\in(0,1)$, and the following compatibility conditions:
\begin{equation}\label{compatibility}
	\begin{cases}
		(\mu+\kappa)\Delta u_0+(\mu+\lambda-\kappa)\nabla\Div u_0+2\kappa\nabla^{\perp}\omega_0-\nabla P(\rho_0)=\rho_0^{\frac{1}{2}}g_1,\\
		\nu\Delta\omega_0+2\kappa\text{curl}~u_0-4\kappa\omega_0=\rho_0^{\frac{1}{2}}g_2,
	\end{cases}
\end{equation}
with $g_1,~g_2\in L^2$. Then in addition to \eqref{bound}$-$\eqref{decay}, the strong solution $(\rho, u, \omega)$ obtained in Theorem \ref{thm1.1} satisfies for any $0<T<\infty$
\begin{align}\label{regularity1}
\left\{
\begin{array}{l}
\nabla^2\rho, \nabla^2P\in C([0,T];L^2\cap L^q),
\\
\bar x^{\delta_0}\nabla^2\rho,\, \bar x^{\delta_0}\nabla^2P, \, \nabla^2{u}, \,\nabla^2{\omega}\in L^\infty(0,T;L^2),
\\
\sqrt{\rho}u_t,\,\sqrt{\rho}\omega_t,\, \sqrt{t}\nabla{u_t},\, \sqrt{t}\nabla{\omega_t}, \, \sqrt{t}{\omega_t}, \, \sqrt{t}\bar{x}^{-1}u_t,\, \sqrt{t}\bar{x}^{-1}\omega_t \in L^\infty(0,T;L^2),
\\
t\sqrt{\rho}u_{tt}, \, t\sqrt{\rho}\omega_{tt}, \, t\nabla^2{u_t},\,
t\nabla^2{\omega_t} \in L^\infty(0,T;L^2),
\\
t\nabla^3u,\, t\nabla^3\omega \in L^\infty(0,T;L^2\cap L^q),
\\
\nabla{u_t}, \, \nabla\omega_t,\, \omega_t, \, \bar{x}^{-1}u_t,\, \bar{x}^{-1}\omega_t, \, t\nabla{u_{tt}},\, t\nabla\omega_{tt},\, t\omega_{tt},\, t\bar{x}^{-1}u_{tt}, t\bar{x}^{-1}\omega_{tt}\in L^2(0,T; L^2),
\\
t\nabla^2(\rho u), t\nabla^2(\rho \omega)\in L^\infty(0,T; L^{\frac{q+2}{2}}).
\end{array}
\right.
\end{align}
\end{theorem}

\begin{remark}
 The compatibility conditions \eqref{compatibility} are necessary for solutions  with the higher regularity stated in \eqref{regularity1}.
 They provide the initial $L^2$ bounds for $\rho^{\frac{1}{2}}\dot u$ and $\rho^{\frac{1}{2}}\dot\omega$ needed to derive the energy estimates for the material derivatives in Lemma \ref{lem_rho_u_t_omega_t_H^1}.
\end{remark}

We briefly describe the main ideas of the proof. Our result concerns the finite mass Cauchy problem with far field vacuum and constant viscosities. In contrast to the results in \cite{Wan_Zhang,Du,Lin_Liu_Zhang}, we assume small initial energy, but do not require density dependent bulk viscosity, large bulk viscosity, or a preparation condition on $\Div u_0$.  Thus, the result addresses a different regime from the large data results cited above. The proof builds on the small energy framework of Li and Xin \cite{Li_Xin}; the additional work is to close the estimates for the coupled velocity and microrotation in the presence of vacuum.

The energy controls $\rho^{\frac{1}{2}}\omega$, whereas the momentum equation contains the unweighted term $2\kappa\nabla^\perp\omega$. We use the combined dissipation to control $\omega$ and the Helmholtz equation \eqref{Helmholtz} to estimate its spatial derivatives. Together with the modified vorticity, these estimates require simultaneous bounds for $\rho\dot u$ and $\rho\dot\omega$; see Lemma \ref{lem_F}. When estimating the material derivatives, the mixed terms are treated in the combined energy estimate, whose dissipation controls $\nabla\dot u$, $\nabla\dot\omega$, and $\dot\omega$. Near the initial time, we also extend the linear decomposition to the pair $(u,\omega)$. In particular, the pressure component $u^P$ forces the angular equation through $\text{curl}~u^P$, so $\omega^P$ must be estimated even though its initial data vanish; see Lemma \ref{lem_u^I_u^P}.

For large times, the coupled estimates yield the decay of the velocity gradient, the pressure, and the microrotational velocity, while the Helmholtz equation gives the decay of $\nabla\omega$ stated in \eqref{decay}. On $\mathbb{R}^2$, we combine these estimates with weighted inequalities and lower bounds for the mass in expanding balls, following \cite{Li_Xin}. The resulting estimates for the effective viscous flux enable us to apply Zlotnik's inequality to obtain a uniform-in-time upper bound for the density and close the small energy argument. Finally, higher regularity is obtained by adapting the estimates of Li and Liang \cite{Li_Liang} to the coupled system.

The rest of this paper is organized as follows. Section \ref{section2} presents some basic inequalities and estimates needed in the subsequent analysis. In Sections \ref{section3} and \ref{section4}, we derive the {\it a priori} estimates required to extend local solutions globally. Finally, in Section \ref{section5}, we give the proofs of Theorems \ref{thm1.1} and \ref{thm1.2}.

\section{Preliminaries}\label{section2}
In this section, we recall some  basic inequalities which will be used frequently below. We begin with a local well-posedness result for strong and classical solutions to problem \eqref{model}$-$\eqref{initial}. This result can be similarly achieved by \cite{Li_Liang}.
\begin{lemma}[Local existence]\label{local}
Assume that the initial data $(\rho_0, u_0, \omega_0)$ satisfy \eqref{initial1} except $u_0\in\dot{H}^\beta$. Then there exists a time $t_0>0$  such that problem \eqref{model}$-$\eqref{initial} admits a unique strong solution $(\rho, u, \omega)$ on $\mathbb{R}^2\times(0,t_0)$ satisfying \eqref{regularity}$-$\eqref{low_bound}. Moreover, if the initial data $(\rho_0, u_0, \omega_0)$ satisfy \eqref{ad_initial} and \eqref{compatibility} besides \eqref{initial1}, then the solution $(\rho, u, \omega)$ also satisfies \eqref{regularity1}.
\end{lemma}

\begin{lemma}[Gagliardo-Nirenberg~\cite{Nirenberg}]\label{Gagliardo}
For $p\in[2,\infty)$, $q\in(1,\infty)$, and $r\in(2,\infty)$, there exists a constant $C>0$ that may depend on $p$, $q$, and $r$  such that for all $f\in H^1$ and $g\in L^q(\mathbb{R}^2)\cap D^{1,r}$, we have
\begin{equation}\label{Gagliardo1}
	\begin{cases}
		\|f\|_{L^p(\mathbb{R}^2)}^p\leq C\|f\|_{L^2(\mathbb{R}^2)}^2\|\nabla f\|_{L^2(\mathbb{R}^2)}^{p-2},
\\[4pt]
\|g\|_{C(\overline{\mathbb{R}^2})}\leq C\|g\|_{L^q(\mathbb{R}^2)}^{\frac{q(r-2)}{2r+q(r-2)}}\|\nabla g\|_{L^r(\mathbb{R}^2)}^{\frac{2r}{2r+q(r-2)}}.
	\end{cases}
\end{equation}
\end{lemma}

The following weighted $L^p$ bounds for functions in $D^1(\mathbb{R}^2)$ can be found in \cite{Lions}.
\begin{lemma}\label{lemma_weight}
 For $m\in[2, \infty)$ and $\theta\in(1+\frac{m}{2},\infty)$, there exists a positive constant $C$ such that for all $v\in D^1({\mathbb{R}^2})$,
\begin{align}\label{weight_L}
\left( \int_{\mathbb{R}^2} \frac{\vert{}v\vert{}^m}{e+\vert{}x\vert{}^2} (\log(e+\vert{}x\vert{}^2))^{-\theta} dx \right)^{1/m} \le C\Vert{}v\Vert{}_{L^2(B_1)} + C\Vert{}\nabla v\Vert{}_{L^2(\mathbb{R}^2)}.
\end{align}
\end{lemma}
Combining Lemma \ref{lemma_weight} with the Poincar\'e inequality yields the following lemma; see \cite{Li_Xin} for a detailed proof.
\begin{lemma}\label{weight_Poincare}
Let $\bar{x}$ be the weight function defined in \eqref{weight}. Assume that $\rho\in L^\infty$ satisfies:
\begin{align}
	0\leq\rho\leq M_1, \quad  M_2\leq\int_{B_{N_*}}\rho dx, \quad \rho\bar{x}^\alpha\in L^1(\mathbb{R}^2),\notag
\end{align}
for $N_*\geq1$ and positive constants $M_1$, $M_2$, and $\alpha$. Then, for $r\in[2,\infty)$, there exists a positive constant $C$, depending only on $M_1$, $M_2$, $\alpha$, and $r$ such that
\begin{align}\label{wp}
\left( \int_{\mathbb{R}^2} \rho \vert{}v\vert{}^r dx \right)^{\frac{1}{r}} \le C N_*^3 (1 + \Vert{}\rho \bar{x}^\alpha\Vert{}_{L^1(\mathbb{R}^2)}) \left( \Vert{}\rho^{\frac{1}{2}} v\Vert{}_{L^2(\mathbb{R}^2)} + \Vert{}\nabla v\Vert{}_{L^2(\mathbb{R}^2)} \right),
\end{align}
for each $v\in\{v\in D^1(\mathbb{R}^2)|\rho^\frac{1}{2}v\in L^2(\mathbb{R}^2)\}$.
\end{lemma}

Next, we derive some basic estimates. First, it follows from the momentum equation in \eqref{model} that the effective viscous flux $F$ and the modified vorticity $\bar\omega$ satisfy the following elliptic equations:
\begin{align}\label{elliptic}
\Delta{F}=\Div(\rho\dot{u}),\quad (\mu+\kappa)\Delta\bar\omega=\text{curl}(\rho\dot{u}).
\end{align}
Furthermore, it follows from the microrotational equation in \eqref{model} that $\omega$ satisfies the following Helmholtz equation:
\begin{align}\label{Helmholtz}
-\nu\Delta\omega+\frac{4\kappa\mu}{\mu+\kappa}\omega=2\kappa\bar\omega-\rho\dot\omega.
\end{align}
\begin{lemma}\label{lem_F}
Let $(\rho, u, \omega)$ be a smooth solution to system \eqref{model}. Then for $2\leq p<\infty$ there exists a positive constant $C$ depending only on $\mu$, $\lambda$, $\kappa$, $\nu$, and $p$ such that
\begin{align}
&\|\nabla{F}\|_{L^p(\mathbb{R}^2)}+\|\nabla{\bar\omega}\|_{L^p(\mathbb{R}^2)}\leq C\|\rho\dot{u}\|_{L^p(\mathbb{R}^2)}, \label{nabla_F_L^p}
\\
&\|{F}\|_{L^p(\mathbb{R}^2)}+\|{\bar\omega}\|_{L^p(\mathbb{R}^2)}\leq C\|\rho\dot{u}\|_{L^2(\mathbb{R}^2)}^{\frac{p-2}{p}}\left(\|\nabla{u}\|_{L^2(\mathbb{R}^2)}+\|{P}\|_{L^2(\mathbb{R}^2)}+\|{\omega}\|_{L^2(\mathbb{R}^2)}\right)^{\frac{2}{p}}, \label{F_L^p}
\\
&\|\nabla\omega\|_{L^p(\mathbb{R}^2)}\leq C\|\nabla\omega\|_{L^2(\mathbb{R}^2)}^{\frac{2}{p}}\left(\|\rho\dot{\omega}\|_{L^2(\mathbb{R}^2)}+\|\rho\dot{u}\|_{L^2(\mathbb{R}^2)}\right)^{\frac{p-2}{p}},  \label{nabla_omega_L^p}
\\
&\|\nabla{u}\|_{L^p(\mathbb{R}^2)}\leq
C\|\rho\dot{u}\|_{L^2(\mathbb{R}^2)}^{\frac{p-2}{p}}
\left(\|\nabla{u}\|_{L^2(\mathbb{R}^2)}+\|{P}\|_{L^2(\mathbb{R}^2)}+\|{\omega}\|_{L^2(\mathbb{R}^2)}\right)^{\frac{2}{p}} \notag
\\
&\qquad\qquad\quad\;\;+C\|{P}\|_{L^p(\mathbb{R}^2)}+C\|\omega\|_{L^p(\mathbb{R}^2)}. \label{nabla_u_L^p}
\end{align}
\end{lemma}
\begin{proof}
On the one hand, the standard $L^p$ estimates for the elliptic system \eqref{elliptic} directly yield \eqref{nabla_F_L^p}. Combining \eqref{nabla_F_L^p} with \eqref{Gagliardo1} and \eqref{flux_curl}, we obtain \eqref{F_L^p}. For $\omega$, multiplying \eqref{Helmholtz} by $-\Delta\omega$ and integrating by parts, we get
\begin{align*}
\nu\|\Delta\omega\|_{L^2}^2+\frac{4\kappa\mu}{\mu+\kappa}\|\nabla\omega\|_{L^2}^2=2\kappa\int\nabla\bar\omega\cdot\nabla\omega dx+\int\rho\dot\omega\Delta\omega dx.
\end{align*}
Together with \eqref{nabla_F_L^p} for $p=2$ and  Young's inequality,  this identity gives
\begin{align}\label{omega_elliptic_L^2}
\|\nabla\omega\|_{L^2}+\|\nabla^2\omega\|_{L^2}\leq C(\|\rho\dot{\omega}\|_{L^2}+\|\rho\dot{u}\|_{L^2}).
\end{align}
Combining \eqref{omega_elliptic_L^2} with \eqref{Gagliardo1}, we obtain \eqref{nabla_omega_L^p}. On the other hand, based on $-\Delta{u}=-\nabla\Div{u}+\nabla^{\perp}\text{curl}~u$, we get
\begin{align}
\nabla{u}=-\nabla(-\Delta)^{-1}\nabla\Div{u}-\nabla(-\Delta)^{-1}(-\nabla^{\perp}\text{curl}~u).\notag
\end{align}
Applying the standard $L^p$ estimates to this identity and using \eqref{F_L^p}, we have
\begin{equation}
	\begin{aligned}
		[b]
		\|\nabla{u}\|_{L^p}&\leq C(p)(\|\Div{u}\|_{L^p}+\|\text{curl}~u\|_{L^p})\\
		&\leq C(p)(\|F\|_{L^p}+\|\bar\omega\|_{L^p}+\|P\|_{L^p}+\|\omega\|_{L^p})\\
		&\leq C(p)\|\rho\dot{u}\|_{L^2}^{\frac{p-2}{p}}\left(\|\nabla{u}\|_{L^2}+\|{P}\|_{L^2(\mathbb{R}^2)}+\|{\omega}\|_{L^2}\right)^{\frac{2}{p}}+C(p)\|{P}\|_{L^p}+C(p)\|\omega\|_{L^p}.\notag
	\end{aligned}
\end{equation}
This completes the proof of the lemma.
\end{proof}

In order to obtain a uniform-in-time upper bound for $\rho$, we need the following Zlotnik inequality.
\begin{lemma}[\cite{Zlotnik}]\label{Zlotnik}
Assume that a function $y(t)$ on $t\in [0, T]$ satisfies
\begin{align}
y'(t)=g(y)+b'(t), \quad y(0)=y^0,\notag
\end{align}
with $g\in C(\mathbb{R})$ and $y,b\in W^{1,1}(0, T)$. If $g(\infty)=-\infty$ and
\begin{align}
b(t_2)-b(t_1)\le N_0+N_1(t_2-t_1),\notag
\end{align}
for all $0\le t_1<t_2\le T$ with some $N_0\ge0$ and $N_1\ge0$, then
\begin{align}
y(t)\le\max\left\{y^0, \bar{\zeta}\right\}+N_0<\infty, ~~~~ t\in [0, T],\notag
\end{align}
where $\bar{\zeta}$ is a constant such that
\begin{align}
g(\zeta)\le-N_1 \quad \text{for all} \quad \zeta\ge\bar{\zeta}.\notag
\end{align}
\end{lemma}
To estimate $\|\nabla u\|_{L^\infty}$ and ~$\|\nabla\rho\|_{L^2\cap L^q}$ for $(q>2)$, we need the following Beale-Kato-Majda type inequality.
\begin{lemma}[\cite{Beale_Kato_Majda, Li_Xin}]\label{lem_BKM}
	For $2<q<\infty$, there exists a positive constant $C(q)$ such that the following estimate holds for all $\nabla u\in L^2(\mathbb{R}^2)\cap D^{1,q}(\mathbb{R}^2)$,
	\begin{equation}\label{BKM}
		\begin{aligned}
			[b]
			\Vert{}\nabla u\Vert{}_{L^\infty(\mathbb{R}^2)}\le& C\left(\Vert{}\Div{u}\Vert{}_{L^\infty(\mathbb{R}^2)}+\Vert{}\operatorname{curl}u\Vert{}_{L^\infty(\mathbb{R}^2)}\right)\log(e+\Vert{}\nabla^2u\Vert{}_{L^q(\mathbb{R}^2)})
\\
&+C\Vert{}\nabla u\Vert{}_{L^2(\mathbb{R}^2)}+C.
		\end{aligned}\notag
	\end{equation}
\end{lemma}

\section{Lower order {\it a priori} estimates}\label{section3}
In this section, we derive the lower-order {\it a priori} estimates for smooth solutions to the Cauchy problem \eqref{model}$-$\eqref{initial} and establish a uniform-in-time upper bound for the density. These estimates play a crucial role in extending the local strong and classical solutions obtained in Lemma \ref{local}. Let $T>0$ be a fixed time and $(\rho, u, \omega)$ be a smooth solution to \eqref{model}$-$\eqref{initial} on $\mathbb{R}^2\times(0,T]$ with smooth initial data
$(\rho_0, u_0, \omega_0)$ satisfying \eqref{initial1}$-$\eqref{initial2}. We first derive the following key {\it a priori} estimates.
\begin{proposition}\label{prop3.1}
Let ~$\sigma(t)\triangleq\min\{1,t\}$. Under the assumptions of Theorem
\ref{thm1.1}, there exists a positive constant $\varepsilon$ depending  on $\mu$, $\lambda$, $\nu$, $\kappa$, $\gamma$, $a$, $\bar\rho$, $\beta$, $N_0$, and $M$ such that if $(\rho, u, \omega)$ is a smooth solution of \eqref{model}$-$\eqref{initial} on $\mathbb{R}^2\times(0,T]$ satisfying
\begin{align}\label{priori_condition}
	\sup_{\mathbb{R}^2\times(0,T]}\rho\leq2\bar\rho, \quad A_1(T)+A_2(T)\leq2E_0^{\frac{1}{2}},
\end{align}
then, provided that $E_0\leq\varepsilon$, the following estimates hold
\begin{align}\label{uniform_bound}
\sup_{\mathbb{R}^2\times(0,T]}\rho\leq\frac{7}{4}\bar\rho, \quad A_1(T)+A_2(T)+\int_0^T\sigma\|P\|_{L^2}^2dt\leq E_0^{\frac{1}{2}},
\end{align}
where
\begin{align}\label{A_1}
A_1(T)\triangleq\sup_{0\leq t\leq T}\left(\sigma(\|\nabla{u}\|_{L^2}^2+\|\nabla{\omega}\|_{L^2}^2+\|\omega\|_{L^2}^2)\right)+\int_0^T\sigma\int\left(\rho|\dot{u}|^2+\rho|\dot{\omega}|^2\right)dxdt,
\end{align}
\begin{align}\label{A_2}
	A_2(T)\triangleq\sup_{0\leq t\leq T}\sigma^2\int\left(\rho|\dot{u}|^2+\rho|\dot{\omega}|^2\right)dx+\int_0^T\sigma^2\int\left(|\nabla\dot{u}|^2+|\nabla\dot{\omega}|^2+|\dot{\omega}|^2\right)dxdt.
\end{align}
\end{proposition}
The complete proof of Proposition \ref{prop3.1} follows from Lemmas \ref{lem_energy_estimate}$-$\ref{lem_rho_infty_bound} below. For convenience, throughout the rest of this paper, we use $C$ to denote a generic positive constant depending on $\mu$, $\lambda$, $\nu$, $\kappa$, $\gamma$, $a$, $\bar\rho$, $\beta$, $N_0$, and $M$, and use $C(\alpha)$ for such a constant depending additionally on $\alpha$.

\begin{lemma}\label{lem_energy_estimate}
Let $(\rho, u, \omega)$ be a smooth solution to \eqref{model}$-$\eqref{initial} on $\mathbb{R}^2\times(0,T]$. Then there exists a positive constant $C$, depending only on $\mu$, $\lambda$, $\kappa$, $\nu$, $\gamma$ such that
\begin{align}
	&\sup_{0\leq t\leq T}\int\left(\rho|u|^2+\rho|\omega|^2+P\right)dx+\int_0^T\int\left(|\nabla{u}|^2+|\nabla{\omega}|^2+|\omega|^2\right)dxdt\leq CE_0, \label{energy_estimate}
	\\
	&A_1(T)\leq C(\bar\rho)E_0+C(\bar\rho)\int_0^T\sigma\left(\|\nabla{u}\|_{L^3}^3+\|\nabla{\omega}\|_{L^3}^3+\|{\omega}\|_{L^3}^3\right)dt. \label{A_1_estimate}
\end{align}
\end{lemma}
\begin{proof}
Multiplying \eqref{model}$_2$ and \eqref{model}$_3$ by $u$ and $\omega$, respectively, integrating by parts over $\mathbb{R}^2$, and adding the resulting identities, we have
\begin{equation}
	\label{energy_estimate1}
	\begin{aligned}
		[b]
		&\frac{d}{dt}\int\left(\frac{1}{2}\rho|u|^2+\frac{1}{2}\rho|\omega|^2+\frac{1}{\gamma-1}P\right)dx+(2\mu+\lambda)\int|\Div u|^2dx
\\
&+(\mu+\kappa)\int|\bar\omega|^2dx
+\frac{4\kappa\mu}{\mu+\kappa}\int|\omega|^2dx+\nu\int|\nabla\omega|^2dx=0,
	\end{aligned}
\end{equation}
	where we have used the following facts
\begin{equation}\label{Hodge}
	\begin{cases}
		\displaystyle\int\bigg[(\mu+\kappa)|\nabla{u}|^2+(\mu+\lambda-\kappa)|\Div u|^2\bigg]dx=\int\bigg[(2\mu+\lambda)|\Div u|^2+(\mu+\kappa)|\text{curl}~u|^2\bigg]dx,
\\[10pt]
\displaystyle\int\bigg[(\mu+\kappa)|\text{curl}~u|^2-4\kappa\omega\text{curl}~u+4\kappa|\omega|^2\bigg]dx
=\int\left[(\mu+\kappa)|\bar\omega|^2+\frac{4\kappa\mu}{\mu+\kappa}|\omega|^2\right]dx.
	\end{cases}
\end{equation}
Further, using
\begin{align*}
\begin{split}
\|\nabla{u}\|_{L^2}^2=\|\Div{u}\|_{L^2}^2+\|\text{curl}~u\|_{L^2}^2,
\quad \|\text{curl}~u\|_{L^2}^2\leq C\left(\|\bar\omega\|_{L^2}^2+\|\omega\|_{L^2}^2\right),
\end{split}
\end{align*}
and integrating \eqref{energy_estimate1} over $[0,T]$, we obtain \eqref{energy_estimate}.

Next, multiplying \eqref{model}$_2$ and \eqref{model}$_3$ by $\dot{u}$ and $\dot{\omega}$, respectively, summing the resulting equations and integrating over  $\mathbb{R}^2$, we derive
\begin{equation}\label{sum}
	\begin{aligned}
		[b]
		&\int\left(\rho|\dot{u}|^2+\rho|\dot{\omega}|^2\right)dx\\
		=&-\int\dot{u}\cdot\nabla{P}dx+(\mu+\kappa)\int\dot{u}\cdot\Delta{u}dx+(\mu+\lambda-\kappa)\int\dot{u}\cdot\nabla\Div{u}dx
	\\
	&+\nu\int\dot{\omega}\Delta\omega dx+2\kappa\int\dot{\omega}\text{curl}~udx-4\kappa\int\omega\dot{\omega}dx+2\kappa\int\dot{u}\cdot\nabla^\perp\omega dx
	\\
	\triangleq&\sum_{i=1}^{7}I_i.
	\end{aligned}
\end{equation}
We now estimate $I_i$. First, by calculations similar to those for $i=1,2,3$ in the proof of Lemma $3.2$ in \cite{Li_Xin}, we get
\begin{equation}
	\label{sum1}
	\begin{aligned}
		[b]
		\sum_{i=1}^{3}I_i\leq&\frac{d}{dt}\left(\int P\Div{u}dx-\frac{\mu+\kappa}{2}\int|\nabla{u}|^2dx-\frac{\mu+\lambda-\kappa}{2}\int|\Div{u}|^2dx\right)
\\
&+C\left(\int P|\nabla{u}|^2dx+\int|\nabla{u}|^3dx\right).
	\end{aligned}
\end{equation}
Then, integrating $I_4$ by parts and combining it with $I_i~(i=5,6,7)$, we have
\begin{equation}\label{sum2}
	\begin{aligned}
		[b]
		\sum_{i=4}^{7}I_i\leq&\frac{d}{dt}\left(-\frac{\nu}{2}\int|\nabla\omega|^2dx+2\kappa\int\omega\text{curl}~udx-2\kappa\int|\omega|^2dx\right)
\\
&+C\int(|\omega|^2+|\nabla\omega|^2)|\nabla{u}|dx.
	\end{aligned}
\end{equation}
Substituting \eqref{sum1} and \eqref{sum2} into \eqref{sum}, and utilizing \eqref{Hodge} together with H\"older inequality, we deduce that
\begin{equation}\label{BD}
	\begin{aligned}[b]
&B'(t)+\int\left(\rho|\dot{u}|^2+\rho|\dot{\omega}|^2\right)dx
\\
\leq& C\left(\int P|\nabla{u}|^2dx+\int|\nabla{u}|^3dx+\int|\nabla{u}||\nabla\omega|^2dx+\int|\omega|^2|\nabla{u}|dx\right),
	\end{aligned}
\end{equation}
where
\begin{equation}\label{B}
	\begin{aligned}[b]
B(t)=&\frac{2\mu+\lambda}{2}\|\Div{u}\|_{L^2}^2+\frac{\mu+\kappa}{2}\|\text{curl}~u\|_{L^2}^2+\frac{\nu}{2}\|\nabla{\omega}\|_{L^2}^2+2\kappa\|{\omega}\|_{L^2}^2
\\
&-2\kappa\int\omega\text{curl}~udx-\int P\Div{u}dx.
	\end{aligned}\notag
\end{equation}
Obviously, $B(t)$ satisfies
\begin{align}\label{B_up}
B(t)\leq C\left(\|\nabla{u}\|_{L^2}^2+\|\nabla{\omega}\|_{L^2}^2+\|{\omega}\|_{L^2}^2+\|{P}\|_{L^2}^2\right).
	\end{align}
In addition, by Young's inequality, we have
\begin{align*}
&\left|\int P\Div{u}dx\right|\leq\frac{\mu+\kappa}{8}\|\Div{u}\|_{L^2}^2+C\|{P}\|_{L^2}^2,
\\
 &\left|2\kappa\int\omega\text{curl}~udx\right|\leq\frac{\mu+2\kappa}{4}\|\text{curl}~{u}\|_{L^2}^2+\frac{4\kappa^2}{\mu+2\kappa}\|{\omega}\|_{L^2}^2.
\end{align*}
Since the viscosity condition gives $2\mu+\lambda\geq\mu+\kappa$, we derive the lower bound of $B(t)$:
\begin{align}\label{B_low}
	\frac{\mu}{4}\|\nabla{u}\|_{L^2}^2+\frac{\nu}{2}\|\nabla{\omega}\|_{L^2}^2
+\frac{2\kappa\mu}{\mu+2\kappa}\|{\omega}\|_{L^2}^2-C\|{P}\|_{L^2}^2\leq B(t).
\end{align}
 Then, multiplying \eqref{BD} by $\sigma$, integrating over $(0,T)$, and employing \eqref{B_up}, \eqref{B_low} and \eqref{energy_estimate}, we have
\begin{align*}
	\begin{split}
A_1(T)\leq& CE_0+C\sup_{0\leq t\leq T}\|P\|_{L^2}^2
\\
&+C\int_0^T\sigma\int\left(P|\nabla{u}|^2+|\nabla{u}|^3+|\nabla{u}||\nabla{\omega}|^2+|\nabla{u}||\omega|^2\right)dxdt.
	\end{split}
\end{align*}
Finally, combining this estimate with \eqref{priori_condition}, \eqref{energy_estimate}, and H\"{o}lder inequality yields \eqref{A_1_estimate}.
\end{proof}

\begin{lemma}\label{lem_A_2_estimate}
Let $(\rho, u, \omega)$ be a smooth solution to \eqref{model}$-$\eqref{initial} on $\mathbb{R}^2\times(0,T]$. Then there exists a constant $C$, depending only on $\mu$, $\lambda$, $\kappa$, $\nu$, and $\gamma$ such that
\begin{align}\label{A_2_estimate}
	\begin{split}
A_2&\leq CA_1(T)+CA_1(T)E_0+C\int_0^T\sigma^2\left(\|\nabla{u}\|_{L^4}^4+\|\nabla{\omega}\|_{L^4}^4+\|{\omega}\|_{L^4}^4+\|{P}\|_{L^4}^4\right)dt.
	\end{split}
\end{align}
\end{lemma}
\begin{proof}
First, applying $\partial_t+\Div(u\cdot{})$ to \eqref{model}$_2^j$,  multiplying the result by $\dot{u}_j$,  summing with respect to $j$, integrating over $\mathbb{R}^2$, and using integration by parts and Young's inequality, we obtain
\begin{equation}\label{A_2_u_part}
	\begin{aligned}[b]
&\frac{1}{2}\frac{d}{dt}\int\rho|\dot{u}|^2dx+(\mu+\kappa)\int|\nabla\dot{u}|^2dx+(\mu+\lambda-\kappa)\int|\Div\dot{u}|^2dx-2\kappa\int\dot\omega\text{curl}~\dot{u}dx
\\
=&(\mu+\kappa)\int\dot{u}_j\partial_k\left(-\partial_ku\cdot\nabla{u}_j+\Div{u}\partial_ku_j\right)dx-(\mu+\kappa)\int\dot{u}_j\Div{(\partial_ku\partial_ku_j)}dx
\\
&-(\mu+\lambda-\kappa)\int\dot{u}_j\partial_j\left(\partial_iu\cdot\nabla{u}_i-(\Div{u})^2\right)dx-(\mu+\lambda-\kappa)\int\dot{u}_j\Div{(\partial_ju\Div{u})}dx
\\
&+(\gamma-1)\int\dot{u}_j\partial_j(P\Div{u})dx+\int\dot{u}_j\Div(P\partial_j{u})dx+2\kappa\int\dot{u}_j(\nabla^\perp\omega\cdot\nabla)u_jdx
\\
\leq& C\|\nabla\dot{u}\|_{L^2}\|\nabla{u}\|_{L^4}^2+C\|\nabla\dot{u}\|_{L^2}\|P\|_{L^4}\|\nabla{u}\|_{L^4}
+2\kappa\int \dot{u}_j\nabla^\perp\omega\cdot\nabla{u}_jdx
\\
\leq&\delta\|\nabla\dot{u}\|_{L^2}^2+C(\delta)(\|\nabla{u}\|_{L^4}^4+\|P\|_{L^4}^4+\|\omega\|_{L^4}^4).
	\end{aligned}
\end{equation}
In the last term of \eqref{A_2_u_part}, we integrate by parts to move the derivative of $\omega$.  The terms containing second derivatives of
$u_j$ cancel because mixed derivatives commute. The remaining factor $\nabla\dot u$ can be absorbed in the dissipation. Indeed,
\begin{equation}\label{new_estimate}
	\begin{aligned}[b]
\Big|2\kappa\int\dot{u}_j\nabla^\perp\omega\cdot\nabla{u}_jdx\Big|=&2\kappa\Big|\int\dot{u}_j(\partial_2\omega\partial_1u_j-\partial_1\omega\partial_2u_j)dx\Big|
\\
=&2\kappa\Big|\int\omega(\partial_1\dot{u}_j\partial_2 u_j-\partial_2\dot{u}_j\partial_1 u_j)dx\Big|
\\
\leq&2\kappa\|\nabla\dot{u}\|_{L^2}\|\nabla u\|_{L^4}\|\omega\|_{L^4}\\
\leq&\delta\|\nabla\dot{u}\|_{L^2}^2+C(\delta)(\|\nabla u\|_{L^4}^4+\|\omega\|_{L^4}^4).
	\end{aligned}\notag
\end{equation}
Next, we use calculations analogous to those for \eqref{A_2_u_part} to obtain from $\eqref{model}_3$ that
\begin{equation}\label{A_2_omega_part}
\begin{aligned}[b]
&\frac{1}{2}\frac{d}{dt}\int\rho|\dot{\omega}|^2dx+4\kappa\int|\dot\omega|^2dx+\nu\int|\nabla{\dot{\omega}}|^2dx-2\kappa\int\dot{\omega}\text{curl}~\dot{u}dx
\\
=&-4\kappa\int\dot{\omega}\omega\Div{u}dx+\nu\int\dot{\omega}\partial_k(-\partial_ku\cdot\nabla\omega+\partial_k\omega\Div{u})dx-\nu\int\dot{\omega}\Div(\partial_ku\partial_k\omega)dx
\\
\leq& C\|\nabla{u}\|_{L^4}\|\omega\|_{L^4}\|\dot\omega\|_{L^2}+C\|\nabla\dot\omega\|_{L^2}\|\nabla{u}\|_{L^4}\|\nabla{\omega}\|_{L^4}
\\
\leq&\delta(\|\dot\omega\|_{L^2}^2+\|\nabla\dot\omega\|_{L^2}^2)+C(\delta)(\|\nabla{u}\|_{L^4}^4+\|\omega\|_{L^4}^4+\|\nabla{\omega}\|_{L^4}^4).
\end{aligned}
\end{equation}
Adding \eqref{A_2_u_part} and \eqref{A_2_omega_part}, choosing $\delta$ sufficiently small, and proceeding as in \eqref{Hodge}, we deduce that
\begin{equation}\label{A_2_inequality}
	\begin{aligned}[b]
&\frac{d}{dt}\left(\int\left(\rho|\dot{u}|^2+\rho|\dot{\omega}|^2 \right)dx\right)+\left(\|\nabla\dot{u}\|_{L^2}^2+\|\nabla\dot{\omega}\|_{L^2}^2+\|\dot{\omega}\|_{L^2}^2\right)
\\
\leq& C\left(\|\nabla{u}\|_{L^4}^4+\|\omega\|_{L^4}^4+\|\nabla{\omega}\|_{L^4}^4+ \|P\|_{L^4}^4\right).
	\end{aligned}
\end{equation}
Multiplying this inequality by $\sigma^2$, integrating over $(0,T)$, and using \eqref{A_1}, we arrive at \eqref{A_2_estimate}.
\end{proof}
The following lemma plays a crucial role in establishing the decay estimates.
\begin{lemma}\label{lem_A_1_and_A_2}
Let $(\rho, u, \omega)$ be a smooth solution to system \eqref{model}$-$\eqref{initial} on $\mathbb{R}^2\times(0,T]$, satisfying ~\eqref{priori_condition}. Then there exists a positive constant $C$ depending only on $\mu$, $\lambda$, $\kappa$, $\nu$, $\gamma$ and $\bar\rho$, such that
\begin{align}\label{A_1_and_A_2}
A_1(T)+A_2(T)+\int_0^T\sigma\|P\|_{L^2}^2dt\leq C(\bar\rho)E_0.
\end{align}
\end{lemma}
\begin{proof}
By virtue of \eqref{Gagliardo1}, \eqref{nabla_omega_L^p}, and \eqref{nabla_u_L^p}, together with \eqref{priori_condition}, \eqref{A_1}, and \eqref{energy_estimate}, we have
	\begin{equation}\label{L^4_estimate}
		\begin{aligned}[b]
&\int_0^T\sigma^2\left(\|\nabla{u}\|_{L^4}^4+\|\nabla{\omega}\|_{L^4}^4+\|P\|_{L^4}^4+\|\omega\|_{L^4}^4\right)dt
\\
\leq& C\int_{0}^{T}\sigma^2\|\rho\dot{u}\|_{L^2}^2\left(\|\nabla{u}\|_{L^2}^2+\|P\|_{L^2}^2+\|\omega\|_{L^2}^2\right)dt+C\int_0^T\sigma^2\|P\|_{L^4}^4dt
\\
&+C\int_0^T\sigma^2\|\nabla\omega\|_{L^2}^2\|\omega\|_{L^2}^2dt+C\int_0^T\sigma^2\|\nabla\omega\|_{L^2}^2\left(\|\rho\dot\omega\|_{L^2}^2+\|\rho\dot u\|_{L^2}^2  \right)dt
\\
\leq& C(\bar\rho)(A_1(T)+E_0)A_1(T)+C(\bar\rho)\int_0^T\sigma^2\|P\|_{L^2}^2dt.
		\end{aligned}
	\end{equation}
To estimate the last term in \eqref{L^4_estimate}, we use \eqref{model}$_2$ to obtain
\begin{align}
P=(-\Delta)^{-1}\Div(\rho\dot{u})+(2\mu+\lambda)\Div{u}.\notag
\end{align}
Then, by H\"{o}lder inequality and the Sobolev inequality, we directly derive
\begin{equation}
	\begin{aligned}[b]
	\|P\|_{L^2}^2&\leq C\|P\|_{L^\frac{4\gamma}{4\gamma-1}}\|(-\Delta)^{-1}\Div{(\rho\dot{u})}\|_{L^{4\gamma}}+C\|P\|_{L^2}\|\nabla{u}\|_{L^2}
	\\
	&\leq C\|\rho\|_{L^1}^{\frac{1}{2}}\|\rho\|_{L^{2\gamma}}^{\gamma-\frac{1}{2}}
	\|\rho\dot{u}\|_{L^\frac{4\gamma}{2\gamma+1}}+C\|P\|_{L^2}\|\nabla{u}\|_{L^2}
	\\
&\leq C\|\rho\|_{L^1}^{\frac{1}{2}}\|\rho\|_{L^{2\gamma}}^{\gamma-\frac{1}{2}}\|\rho\|_{L^{2\gamma}}^\frac{1}{2}\|\rho^\frac{1}{2}\dot{u}\|_{L^2}+C\|P\|_{L^2}\|\nabla{u}\|_{L^2}
	\\
	&\leq C\|P\|_{L^2}\|\rho^\frac{1}{2}\dot{u}\|_{L^2}+C\|P\|_{L^2}\|\nabla{u}\|_{L^2}.
	\end{aligned}\notag
\end{equation}
Here, in the last step, we have used
\begin{align}
\int\rho dx=\int\rho_0dx=1.\notag
\end{align}
Thus, one has
\begin{align}\label{P_estimate}
\int_0^T\sigma\|P\|_{L^2}^2dt\leq C\int_0^T\sigma\left(\|\rho^\frac{1}{2}\dot{u}\|_{L^2}^2+\|\nabla{u}\|_{L^2}^2\right)dt\leq C(A_1(T)+E_0).
\end{align}
For $A_1(T)$, it follows from \eqref{A_1_estimate} that
\begin{equation}\label{A_1_2}
	\begin{aligned}[b]
A_1(T)\leq& C(\bar\rho)E_0+C(\bar\rho)\int_0^{\sigma(T)}\sigma\left(\|\nabla{u}\|_{L^3}^3+\|\nabla{\omega}\|_{L^3}^3+\|{\omega}\|_{L^3}^3\right)dt
\\
&+C(\bar\rho)\int_{\sigma(T)}^T\sigma\left(\|\nabla{u}\|_{L^3}^3+\|\nabla{\omega}\|_{L^3}^3+\|{\omega}\|_{L^3}^3\right)dt.
	\end{aligned}
\end{equation}
For the second term on the right-hand side of \eqref{A_1_2}, we apply \eqref{F_L^p}$-$\eqref{nabla_u_L^p}, \eqref{priori_condition}, \eqref{A_1}, and \eqref{energy_estimate} to yield
\begin{equation}\label{L^3_estimate1}
	\begin{aligned}[b]
&\int_0^{\sigma(T)}\sigma\left(\|\nabla{u}\|_{L^3}^3+\|\nabla{\omega}\|_{L^3}^3+\|{\omega}\|_{L^3}^3\right)dt
\\
\leq &CE_0+C\int_0^{\sigma(T)}\sigma\|\rho\dot{u}\|_{L^2}\left(\|\nabla{u}\|_{L^2}^2+\|P\|_{L^2}^2+\|\omega\|_{L^2}^2\right)dt
\\
&+C\int_0^{\sigma(T)}\sigma\left(\|\rho\dot{\omega}\|_{L^2}+
\|\rho\dot{u}\|_{L^2}\right)\|\nabla\omega\|_{L^2}^2dt
+C\int_0^{\sigma(T)}\sigma\|\omega\|_{L^2}^2\|\nabla\omega\|_{L^2}dt
\\
\leq& CE_0.
	\end{aligned}
\end{equation}
Additionally, applying \eqref{L^4_estimate}, \eqref{P_estimate}, \eqref{priori_condition}, \eqref{energy_estimate} to the third term, one gets
\begin{equation}\label{L^3_estimate2}
\begin{aligned}[b]
&\int_{\sigma(T)}^T\sigma\left(\|\nabla{u}\|_{L^3}^3+\|\nabla{\omega}\|_{L^3}^3+\|{\omega}\|_{L^3}^3\right)dt
\\
\leq& C(\delta)E_0+\delta\int_{\sigma(T)}^T\left(\|\nabla{u}\|_{L^4}^4+\|\nabla{\omega}\|_{L^4}^4+\|\omega\|_{L^4}^4\right)dt
\\
\leq& C(\delta)E_0+\delta C(\bar\rho)E_0+\delta C(\bar\rho)A_1(T).
\end{aligned}
\end{equation}
We now choose $\delta$ small enough and substitute \eqref{L^3_estimate1} and \eqref{L^3_estimate2} into \eqref{A_1_2}. It follows that
\begin{align}\label{uniform_A_1}
A_1\leq C(\bar\rho)E_0.
\end{align}
Eventually, combining \eqref{L^4_estimate}, \eqref{P_estimate}, \eqref{uniform_A_1}, and \eqref{A_2_estimate}, we arrive at \eqref{A_1_and_A_2}.
\end{proof}

In the next lemma, we derive decay rates of $\nabla{u}$, $P$, and $\omega$. These estimates play a critical role in establishing a uniform-in-time upper bound for the density.

\begin{lemma}\label{lem_decay}
For $p\in[2,\infty)$, there exists a positive constant $C(p,\bar\rho)$, depending only on $p$, $\mu$, $\lambda$, $\kappa$, $\nu$, $\gamma$ and $\bar\rho$ such that if $(\rho, u, \omega)$ is a smooth solution to \eqref{model}$-$\eqref{initial} on $\mathbb{R}^2\times(0,T]$, satisfying \eqref{priori_condition}, then
	\begin{align}\label{decay_estimate1}
		\sup_{\sigma(T)\leq t\leq T}\left(t^{p-1}(\|\nabla{u}\|_{L^p}^p+\|P\|_{L^p}^p+\|\omega\|_{L^p}^p)+t^2(\|\rho^{\frac{1}{2}}\dot{u}\|_{L^2}^2+\|\rho^{\frac{1}{2}}\dot{\omega}\|_{L^2}^2)\right)\leq C(p,\bar\rho)E_0,
	\end{align}
and
\begin{align}\label{decay_estimate2}
\sup_{\sigma(T)\leq t\leq T}t^2\left(\|\nabla{\omega}\|_{L^2}^2+\|\nabla^2{\omega}\|_{L^2}^2\right)\leq C(p,\bar\rho)E_0.
\end{align}
\end{lemma}
\begin{proof}
The pressure $P$ satisfies the equation
$$
P_t+u\cdot\nabla{P}+\gamma P\Div{u}=0.
$$
Multiplying this by $pP^{p-1}$ and integrating over $\mathbb{R}^2$, utilizing $\Div{u}=\frac{P+F}{2\mu+\lambda}$, one obtains
\begin{equation}
	\begin{aligned}[b]
\left(\|P\|_{L^p}^p\right)_t+\frac{p\gamma-1}{2\mu+\lambda}\|P\|_{L^{p+1}}^{p+1}&=-\frac{p\gamma-1}{2\mu+\lambda}\int P^pFdx
\\
&\leq\frac{p\gamma-1}{2(2\mu+\lambda)}\|P\|_{L^{p+1}}^{p+1}+C(p)\|F\|_{L^{p+1}}^{p+1},
	\end{aligned}\notag
\end{equation}
which together with \eqref{F_L^p} gives
\begin{align}\label{P_decay_p}
\begin{split}
\left(\|P\|_{L^p}^p\right)_t+\frac{p\gamma-1}{2(2\mu+\lambda)}\|P\|_{L^{p+1}}^{p+1}\leq C(p)\|\rho\dot{u}\|_{L^2}^{p-1}\left(\|\nabla{u}\|_{L^2}^2+\|{P}\|_{L^2}^2+\|{\omega}\|_{L^2}^2\right).
\end{split}
\end{align}
Particularly, taking $p=2$, we have
\begin{align}\label{P_decay}
\left(\|P\|_{L^2}^2\right)_t+\frac{2\gamma-1}{2(2\mu+\lambda)}\|P\|_{L^{3}}^{3}\leq \delta\|\rho\dot{u}\|_{L^2}^2+C(\delta)\left(\|\nabla{u}\|_{L^2}^4+\|{P}\|_{L^2}^4+\|{\omega}\|_{L^2}^4\right).
\end{align}
One deduces from \eqref{BD} and \eqref{nabla_u_L^p} that
\begin{equation}\label{BD1}
	\begin{aligned}[b]
&B'(t)+\int\left(\rho|\dot{u}|^2+\rho|\dot{\omega}|^2\right)dx\leq C(\|\nabla{u}\|_{L^3}^3+\|\nabla{\omega}\|_{L^3}^3+\|{P}\|_{L^3}^3+\|{\omega}\|_{L^3}^3)
\\
\leq& C\|{P}\|_{L^3}^3+C(\bar\rho)\left(\int\Big(\rho|\dot{u}|^2+\rho|\dot{\omega}|^2\Big)dx\right)^{\frac{1}{2}}\left(\|\nabla{u}\|_{L^2}^2+\|\nabla{\omega}\|_{L^2}^2+\|{P}\|_{L^2}^2+\|{\omega}\|_{L^2}^2\right)
\\
\leq& C\|{P}\|_{L^3}^3+\delta\int\Big(\rho|\dot{u}|^2+\rho|\dot{\omega}|^2\Big)dx+C(\bar\rho,\delta)\left(\|\nabla{u}\|_{L^2}^4+\|\nabla{\omega}\|_{L^2}^4+\|{P}\|_{L^2}^4+\|{\omega}\|_{L^2}^4\right).
	\end{aligned}
\end{equation}
Recalling the properties of $B(t)$ in \eqref{B_up} and \eqref{B_low}, we choose a sufficiently large constant $C_1$ such that
\begin{equation}\label{new_B_bound}
	\begin{aligned}[b]
&\frac{\mu}{4}\|\nabla{u}\|_{L^2}^2+\frac{\nu}{2}\|\nabla{\omega}\|_{L^2}^2+\frac{2\kappa\mu}{\mu+2\kappa}\|{\omega}\|_{L^2}^2+\|{P}\|_{L^2}^2\leq B(t)+C_1\|{P}\|_{L^2}^2
\\
\leq& C\left(\|\nabla{u}\|_{L^2}^2+\|\nabla{\omega}\|_{L^2}^2+\|{\omega}\|_{L^2}^2+\|{P}\|_{L^2}^2\right).
\end{aligned}
\end{equation}
Multiplying \eqref{P_decay} by $C_1$, adding it to \eqref{BD1}, and choosing $\delta$ suitably small, we obtain
\begin{equation}\label{B_estimate}
	\begin{aligned}[b]
&\Big(B(t)+C_1\|{P}\|_{L^2}^2\Big)'+\int\left(\rho|\dot{u}|^2+\rho|\dot{\omega}|^2+P^3\right)dx
\\
\leq &C\left(\|\nabla{u}\|_{L^2}^4+\|\nabla{\omega}\|_{L^2}^4+\|{P}\|_
{L^2}^4+\|{\omega}\|_{L^2}^4\right).
	\end{aligned}
\end{equation}
Multiplying \eqref{B_estimate} by $t$ and using Gronwall inequality together with  \eqref{new_B_bound}, \eqref{A_1_and_A_2}, \eqref{energy_estimate}, and \eqref{priori_condition}, we have
\begin{equation}\label{decay_part1}
	\begin{aligned}[b]
&\sup_{\sigma(T)\leq t\leq T}t\Big(\|\nabla{u}\|_{L^2}^2+\|\nabla{\omega}\|_{L^2}^2+\|{P}\|_{L^2}^2+\|{\omega}\|_{L^2}^2\Big)
\\
&+\int_{\sigma(T)}^{T}t\int\left(\rho|\dot{u}|^2+\rho|\dot{\omega}|^2+P^3\right)dxdt
\leq C(\bar\rho)E_0.
	\end{aligned}
\end{equation}
Next, multiplying \eqref{A_2_inequality} by $t^2$, we obtain
\begin{equation}\label{A_2_decay}
	\begin{aligned}[b]
&\bigg(t^2\int\Big(\rho|\dot{u}|^2+\rho|\dot{\omega}|^2\Big)dx\bigg)_t+t^2\left(\|\nabla\dot{u}\|_{L^2}^2+\|\nabla\dot{\omega}\|_{L^2}^2+\|\dot{\omega}\|_{L^2}^2\right)
\\
\leq& C(\bar\rho)t^2\left(\|\rho\dot{u}\|_{L^2}^2+\|\rho\dot{\omega}\|_{L^2}^2\right)\left(\|\nabla{u}\|_{L^2}^2+\|\nabla{\omega}\|_{L^2}^2+\|{P}\|_{L^2}^2+\|{\omega}\|_{L^2}^2\right)
\\
&+2t\int\left(\rho|\dot{u}|^2+\rho|\dot{\omega}|^2\right)dx+C(\bar\rho)t^2\|{P}\|_{L^4}^4.
	\end{aligned}
\end{equation}
Multiplying \eqref{P_decay_p} by $C_2t^2$ with $p=3$ and $C_2$ sufficiently large, and adding the estimate to \eqref{A_2_decay}, we arrive at
\begin{align*}
	\begin{split}
&\bigg(t^2\Big(\int\big(\rho|\dot{u}|^2+\rho|\dot{\omega}|^2\big)dx+C_2\|P\|_{L^3}^3\Big) \bigg)_t+t^2\left(\|\nabla\dot{u}\|_{L^2}^2+\|\nabla\dot{\omega}\|_{L^2}^2+\|\dot{\omega}\|_{L^2}^2+\|P\|_{L^4}^4\right)
\\
\leq &C(\bar\rho)t^2\left(\|\rho^{\frac{1}{2}}\dot{u}\|_{L^2}^2+\|\rho^{\frac{1}{2}}\dot{\omega}\|_{L^2}^2\right)\left(\|\nabla{u}\|_{L^2}^2+\|\nabla{\omega}\|_{L^2}^2+\|{P}\|_{L^2}^2+\|{\omega}\|_{L^2}^2\right)
\\
&+2t\Big(\int\left(\rho|\dot{u}|^2+\rho|\dot{\omega}|^2\right)dx+C_2\|P\|_{L^3}^3\Big),
	\end{split}
\end{align*}
which, along with  Gronwall inequality, \eqref{decay_part1}, and \eqref{priori_condition}, gives
\begin{equation}\label{decay_part2}
	\begin{aligned}[b]
	&\sup_{\sigma(T)\leq t\leq T}t^2\left(\int\left(\rho|\dot{u}|^2+\rho|\dot{\omega}|^2\right)dx+\|P\|_{L^3}^3\right)
	\\
	&+\int_{\sigma(T)}^{T}t^2\left(\|\nabla\dot{u}\|_{L^2}^2+\|\nabla\dot{\omega}\|_{L^2}^2+\|\dot{\omega}\|_{L^2}^2+\|P\|_{L^4}^4 \right)dt\leq C(\bar\rho)E_0.
	\end{aligned}
\end{equation}
Combining \eqref{priori_condition}, \eqref{decay_part2} and \eqref{omega_elliptic_L^2}, one further deduces
\begin{align}\label{omega_decay}
\sup_{\sigma(T)\leq t\leq T}t^2\left(\|\nabla\omega\|_{L^2}^2+\|\nabla^2\omega\|_{L^2}^2\right)\leq C(\bar\rho)E_0,
\end{align}
which is \eqref{decay_estimate2}.

We now prove the following estimate by induction,
\begin{align}\label{induction}
\sup_{\sigma(T)\leq t\leq T}t^m\|P\|_{L^{m+1}}^{m+1}+\int_{\sigma(T)}^Tt^m\|P\|_{L^{m+2}}^{m+2}dt\leq C(m,\bar\rho)E_0.
\end{align}
Indeed, \eqref{decay_part1} shows that \eqref{induction} holds for  $m=1$. Suppose that \eqref{induction} holds for $m=n$. Then we have
\begin{align}\label{induction_n+1}
\sup_{\sigma(T)\leq t\leq T}t^n\|P\|_{L^{n+1}}^{n+1}+\int_{\sigma(T)}^Tt^n\|P\|_{L^{n+2}}^{n+2}dt\leq C(n,\bar\rho)E_0.
\end{align}
Multiplying \eqref{P_decay_p} by $t^{n+1}$ with $p=n+2$, we have
\begin{equation}\label{3.45}
	\begin{aligned}[b]
&\frac{2(2\mu+\lambda)}{(n+2)\gamma-1}\left(t^{n+1}\|P\|_{L^{n+2}}^{n+2}\right)_t+t^{n+1}\|P\|_{L^{n+3}}^{n+3}
\\
\leq &C(n,\bar{\rho})t^n\|P\|_{L^{n+2}}^{n+2}+C(n,\bar\rho)E_0(\|\nabla{u}\|_{L^2}^2+\|P\|_{L^2}^2+\|\omega\|_{L^2}^2),
	\end{aligned}
\end{equation}
where we have used \eqref{decay_part2}.

Integrating \eqref{3.45} over $[\sigma(T),T]$, and combining the result with \eqref{induction_n+1} and \eqref{A_1_and_A_2}, we obtain \eqref{induction} for $m=n+1$.

Thus, the induction argument yields \eqref{induction}. Finally, \eqref{Gagliardo1}, \eqref{nabla_u_L^p}, \eqref{decay_part1}, \eqref{decay_part2}, \eqref{omega_decay}, and \eqref{induction} lead to \eqref{decay_estimate1}. This completes the proof of this lemma.
\end{proof}

\begin{remark}
	The pressure decay estimate extends to $r>1$. Indeed, interpolation between \eqref{energy_estimate} and \eqref{induction} yields
 $$
 \sup_{\sigma(T)\le t\le T}t^{r-1}\|P(t)\|_{L^r}^{r}\le C(r,\bar\rho)E_0,
\qquad r\in(1,\infty).
 $$
\end{remark}

Combining Lemma \ref{lem_rho_int_low_bound} with Lemma \ref{weight_Poincare} gives the $L^p$ estimates for $\rho\dot{u}$, which are essential for establishing a uniform-in-time upper bound for the density for large times.
\begin{lemma}\label{lem_rho_int_low_bound}
Let $(\rho, u, \omega)$ be a smooth solution to system \eqref{model}$-$\eqref{initial} on $\mathbb{R}^2\times(0,T]$, satisfying the assumptions of Theorem \ref{thm1.1} and the condition \eqref{priori_condition}. Then, for any $\alpha>1$, there exists a positive constant $N_1$, depending only on $\alpha$, $N_0$, and $M$, such that, for all $t\in(0,T]$,
\begin{align}\label{rho_int_low_bound}
\int_{B_{N_1(1+t)\log^{\alpha}(e+t)}}\rho(x,t)dx\geq\frac{1}{4}.
\end{align}
Moreover, for all $t\in[0,T]$, we have
\begin{align}\label{mass_time_bound}
\sup_{0\leq s\leq t}\int\rho(x, s)(1+|x|^2)^{\frac{1}{2}}dx\leq C(M)(1+t).
\end{align}
\end{lemma}
\begin{proof}
The proof can be found in Lemma $3.5$ of \cite{Li_Xin}.
\end{proof}

To establish an upper bound for the density near the initial time, we also need the following lemma.
\begin{lemma}\label{lem_u^I_u^P}
Let $(\rho, u, \omega)$ be a smooth solution to system \eqref{model}$-$\eqref{initial} on $\mathbb{R}^2\times(0,T]$, satisfying \eqref{priori_condition} and the assumptions of Theorem \ref{thm1.1}. Then there exists a positive constant $C$ depending only on $\mu$, $\lambda$, $\kappa$, $\nu$, $\gamma$, $a$, $\bar\rho$, $\beta$, $N_0$ and $M$ such that
\begin{align}\label{u^I_u^P}
\sup_{0\leq t\leq\sigma(T)}t^{1-\beta}\left(\|\nabla{u}\|_{L^2}^2+\|\omega\|_{H^1}^2\right)+\int_0^{\sigma(T)}t^{1-\beta}\int\left(\rho|\dot{u}|^2+\rho|\dot{\omega}|^2\right)dxdt\leq C(\bar\rho, M).
\end{align}
\end{lemma}
\begin{proof}
Motivated by \cite{Hoff,Li_Xin}, we regard $\rho$ and $u$ as given coefficients and introduce the following two linear differential operators $L_1$ and $L_2$ defined by
\begin{equation}
	\begin{cases}
		L_1(v,w)=\rho v_t+\rho u\cdot\nabla{v}-(\mu+\kappa)\Delta{v}-(\mu+\lambda-\kappa)\nabla\Div{v}-2\kappa\nabla^\perp{w},
\\
L_2(v,w)=\rho w_t+\rho u\cdot\nabla{w}+4\kappa{w}-\nu\Delta{w}-2\kappa\text{curl}~v,
	\end{cases}\notag
\end{equation}
where $v(x,t)=(v_1, v_2)(x,t)$ is a vector valued function, while $w(x,t)$ is a scalar valued function. Let $(u^I, \omega^I)$ and $(u^P, \omega^P)$ be the solutions to the following problems:
\begin{align}\label{decompose1}
\left\{
	\begin{array}{l}
L_1(u^I, \omega^I)=0,
\\
L_2(u^I, \omega^I)=0,
\\
u^I|_{t=0}=u_0^I,\quad \omega^I|_{t=0}=\omega_0^I,
	\end{array}
\right.
\end{align}
and
\begin{align}\label{decompose2}
\left\{
	\begin{array}{l}
L_1(u^P, \omega^P)=-\nabla{P},
\\
L_2(u^P, \omega^P)=0,
\\
u^P|_{t=0}=0,\quad \omega^P|_{t=0}=0,
	\end{array}
\right.
\end{align}
 respectively. Applying the standard energy method to \eqref{decompose1}, we obtain
\begin{equation}\label{energy_decompose1}
	\begin{aligned}[b]
&\sup_{0\leq t\leq\sigma(T)}\int\left(\rho|u^I|^2+\rho|\omega^I|^2 \right)dx+\int_0^{\sigma(T)}\left(\|\nabla u^I\|_{L^2}^2+\|\nabla{\omega^I}\|_{L^2}^2+\|\omega^I\|_{L^2}^2\right)dt
\\
\leq& C(\bar\rho)\int\left(|u_0^I|^2+|\omega_0^I|^2\right)dx.
	\end{aligned}
\end{equation}
Next, multiplying \eqref{decompose1}$_1$ and \eqref{decompose1}$_2$ by $\dot{u}^I$ and $\dot\omega^I$, respectively,  summing the resulting identities and performing calculations similar to those in \eqref{Hodge}, we have
\begin{equation}\label{energy_decompose12}
\begin{aligned}[b]
&\frac{1}{2}\frac{d}{dt}\bigg((\mu+\kappa)\|\nabla u^I\|_{L^2}^2+(\mu+\lambda-\kappa)\|\Div u^I\|_{L^2}^2+\nu\|\nabla\omega^I\|_{L^2}^2
\\
&+4\kappa\|\omega^I\|_{L^2}^2-4\kappa\int\omega^I\text{curl}~u^Idx\bigg)
+\int\left(\rho|\dot{u}^I|^2+\rho|\dot\omega^I|^2\right)dx
\\
\leq& C(\bar\rho)\|\nabla{u}\|_{L^2}\left(\|\nabla u^I\|_{L^4}^2+\|\omega^I\|_{L^4}^2+\|\nabla{\omega}^I\|_{L^4}^2+\|\omega^I\|_{L^4}\|\nabla{u}^I\|_{L^4} \right)
\\
\leq&\frac{1}{2}\int\left(\rho|\dot{u}^I|^2+\rho|\dot\omega^I|^2\right)dx+C(M,\bar\rho)\|\nabla{u}\|_{L^2}^2\left(\|\nabla u^I\|_{L^2}^2+\|\omega^I\|_{L^2}^2+\|\nabla{\omega}^I\|_{L^2}^2\right).
\end{aligned}
\end{equation}
Here, the last inequality follows from the standard elliptic estimates
\begin{equation}
	\begin{cases}
		\|\nabla^2u^I\|_{L^2}\leq C(\bar\rho)(\|\rho^{\frac{1}{2}}\dot{u}^I\|_{L^2}+\|\rho^{\frac{1}{2}}\dot{\omega}^I\|_{L^2}),
\\
\|\nabla\omega^I\|_{L^2}+\|\nabla^2\omega^I\|_{L^2}\leq C(\bar\rho)(\|\rho^{\frac{1}{2}}\dot{u}^I\|_{L^2}+\|\rho^{\frac{1}{2}}\dot{\omega}^I\|_{L^2}).
	\end{cases}\notag
\end{equation}
Integrating \eqref{energy_decompose12} over $[0,\sigma(T)]$, and using the equivalence between the quadratic form in \eqref{energy_decompose12} and $\|\nabla u^I\|_{L^2}^2+\|\nabla\omega^I\|_{L^2}^2+\|\omega^I\|_{L^2}^2$, \eqref{energy_decompose1}, \eqref{energy_estimate}, and Gronwall inequality, we achieve that
\begin{equation}\label{interpolation1}
	\begin{aligned}[b]
&\sup_{0\leq t\leq\sigma(T)}(\|\nabla u^I\|_{L^2}^2+\|\nabla\omega^I\|_{L^2}^2+\|\omega^I\|_{L^2}^2)
\\
&+\int_0^{\sigma(T)}\int(\rho|\dot{u}^I|^2+\rho|\dot\omega^I|^2)dxdt
\leq C(\bar\rho,M)(\|\nabla u_0^I\|_{L^2}^2+\|\omega_0^I\|_{H^1}^2),
	\end{aligned}
\end{equation}
and
\begin{equation}\label{interpolation2}
	\begin{aligned}[b]
&\sup_{0\leq t\leq\sigma(T)}t(\|\nabla u^I\|_{L^2}^2+\|\nabla\omega^I\|_{L^2}^2+\|\omega^I\|_{L^2}^2)
\\
&+\int_0^{\sigma(T)}t\int\left(\rho|\dot{u}^I|^2+\rho|\dot\omega^I|^2\right)dxdt
\leq C(\bar\rho,M)(\|u_0^I\|_{L^2}^2+\|\omega_0^I\|_{H^1}^2).
	\end{aligned}
\end{equation}
Since the solution operator $(u_0^I, \omega_0^I)\to(u^I, \omega^I)(\cdot, t)$ is linear, by the standard Stein-Weiss interpolation theorem
\cite{Bergh}, one deduces from \eqref{interpolation1} and \eqref{interpolation2} that for any $\theta\in[\beta,1]$,
\begin{equation}\label{u^I_estimate}
	\begin{aligned}[b]
&\sup_{0\leq t\leq\sigma(T)}t^{1-\theta}\Big(\|\nabla u^I\|_{L^2}^2+\|\nabla\omega^I\|_{L^2}^2+\|\omega^I\|_{L^2}^2\Big)
\\
&+\int_0^{\sigma(T)}t^{1-\theta}\int\Big(\rho|\dot{u}^I|^2+\rho|\dot\omega^I|^2\Big)dxdt
\leq C(\bar\rho,M)(\|u_0^I\|_{\dot{H}^\theta}^2+\|\omega_0^I\|_{H^1}^2),
	\end{aligned}
\end{equation}
where the constant $C$ is independent of $\theta$.

To estimate $(u^P, \omega^P)$, we proceed as in Lemma \ref{lem_F} and obtain
\begin{align}
\begin{split}
&\|\nabla((2\mu+\lambda)\Div{u}^P-P)\|_{L^2}+\|\nabla(\text{curl}~u^P-\frac{2\kappa}{\mu+\kappa}\omega^P)\|_{L^2}\leq C\|\rho^\frac{1}{2}\dot{u}^P\|_{L^2},
\\
&\|\nabla\omega^P\|_{L^2}+\|\nabla^2\omega^P\|_{L^2}\leq C(\|\rho^\frac{1}{2}\dot{u}^P\|_{L^2}+\|\rho^\frac{1}{2}\dot{\omega}^P\|_{L^2}),
\end{split}\notag
\end{align}
and
\begin{equation}\label{u^P_L^4}
	\begin{aligned}[b]
\|\nabla{u}^P\|_{L^4}^2
\leq &C\bigg(\|\Div{u}^P\|_{L^4}^2+\|\text{curl}~u^P\|_{L^4}^2\bigg)
\\
\leq& C\bigg(\|(2\mu+\lambda)\Div{u}^P-P\|_{L^4}^2+\|(\text{curl}~u^P-\frac{2\kappa}{\mu+\kappa}\omega^P)\|_{L^4}^2+\|P\|_{L^4}^2+\|\omega^P\|_{L^4}^2\bigg)
\\
\leq& C(\bar\rho)(\|\rho^\frac{1}{2}\dot{u}^P\|_{L^2}+\|\rho^\frac{1}{2}\dot{\omega}^P\|_{L^2})(1+\|\nabla{u}^P\|_{L^2}+\|\omega^P\|_{L^2})+C(\bar\rho).
	\end{aligned}
\end{equation}
Multiplying \eqref{decompose2}$_1$ and \eqref{decompose2}$_2$ by $\dot{u}^P$ and $\dot\omega^P$, respectively, and adding up the results, one gets
\begin{align}
\begin{split}
&\frac{1}{2}\frac{d}{dt}\bigg((\mu+\kappa)\|\nabla{u}^P\|_{L^2}^2+(\mu+\lambda-\kappa)\|\Div{u}^P\|_{L^2}^2+4\kappa\|\omega^P\|_{L^2}^2+\nu\|\nabla{\omega}^P\|_{L^2}^2
\\
& -4\kappa\int\omega^P\text{curl}~u^Pdx-2\int P\Div{u}^Pdx\bigg)+\int\left(\rho|\dot{u}^P|^2+\rho|\dot{\omega}^P|^2\right)dx
\\
\leq& C\|\nabla{u}\|_{L^2}(\|\nabla u^P\|_{L^4}^2+\|\omega^P\|_{L^4}^2+\|\nabla{\omega}^P\|_{L^4}^2+\|\omega^P\|_{L^4}\|\nabla{u}^P\|_{L^4}+\|P\|_{L^4}\|\nabla u^P\|_{L^4})
\\
\leq&\frac{1}{2}\int\left(\rho|\dot{u}^P|^2+\rho|\dot{\omega}^P|^2\right)dx+C(\bar\rho)\|\nabla{u}\|_{L^2}^2\left(\|\nabla u^P\|_{L^2}^2+\|\omega^P\|_{L^2}^2+\|\nabla{\omega}^P\|_{L^2}^2+1\right)+C(\bar\rho).
\end{split}\notag
\end{align}
Here, the last inequality follows from \eqref{u^P_L^4} and $\|P\|_{L^4}\leq C(\bar\rho)$. Furthermore, using an estimate similar to \eqref{B_low}, together with the fact
\begin{align*}
\left|\int P\Div{u}^Pdx\right|\leq\frac{\mu+\kappa}{4}\|\Div{u}^P\|_{L^2}^2+C\|P\|_{L^2}^2\leq\frac{\mu+\kappa}{4}\|\Div{u}^P\|_{L^2}^2+C(\bar\rho),
\end{align*}
 for all $t\in[0,\sigma(T)]$, we have
\begin{align*}
	\begin{split}
&\|\nabla{u}^P\|_{L^2}^2+\|\nabla{\omega}^P\|_{L^2}^2+\|\omega^P\|_{L^2}^2+\int_0^t\int\left(\rho|\dot{u}^P|^2+\rho|\dot{\omega}^P|^2\right)dxdt
\\
\leq& C(\bar\rho)\int_0^t\|\nabla{u}\|_{L^2}^2\left(\|\nabla u^P\|_{L^2}^2+\|\omega^P\|_{L^2}^2+\|\nabla{\omega}^P\|_{L^2}^2+1\right)dt+C(\bar\rho).
	\end{split}
\end{align*}
Consequently, applying Gronwall  inequality, we obtain
\begin{equation}\label{u^P_estimate}
	\begin{aligned}[b]
&\sup_{0\leq t\leq\sigma(T)}\Big(\|\nabla{u}^P\|_{L^2}^2+\|\nabla{\omega}^P\|_{L^2}^2+\|\omega^P\|_{L^2}^2\Big)
\\
&+\int_0^{\sigma(T)}\int\left(\rho|\dot{u}^P|^2+\rho|\dot{\omega}^P|^2\right)dxdt\leq C(\bar\rho,M).
	\end{aligned}
\end{equation}
Taking $(u_0^I, \omega_0^I)=(u_0,\omega_0)$, we have $u^I+u^P=u$ and $\omega^I+\omega^P=\omega$. Applying \eqref{u^I_estimate} with $\theta=\beta$ and using \eqref{u^P_estimate}, we get \eqref{u^I_u^P}. This completes the proof of Lemma \ref{lem_u^I_u^P}.
\end{proof}
For convenience in the subsequent analysis, we establish the following  weighted estimate.
\begin{lemma}
For any $p\in[2,\infty)$ and ~$t>0$, there exists a positive constant $C$ depending only on $\mu$, $\lambda$, $\kappa$, $\nu$, $\bar\rho$, $p$ and $M$ such that, for any $v\in D^1$ satisfying $\rho^{\frac{1}{2}}v\in L^2$, one derives
\begin{align}\label{nabla_F_L^p1}
	\begin{split}
\|\rho{v}\|_{L^p}\leq C(1+t)^5(\|\rho^{\frac{1}{2}}{v}\|_{L^2}+\|\nabla{v}\|_{L^2}).
	\end{split}
\end{align}
\end{lemma}
\begin{proof}
In fact, combining \eqref{wp}, \eqref{rho_int_low_bound}, \eqref{mass_time_bound}, and \eqref{priori_condition} gives \eqref{nabla_F_L^p1}.
\end{proof}
We now turn to establishing a uniform-in-time upper bound for the density, which is crucial for deriving the higher order estimates and extending the local classical solution globally.
\begin{lemma}\label{lem_rho_infty_bound}
Let $(\rho, u, \omega)$ be a smooth solution to the problem \eqref{model}$-$\eqref{initial} on $\mathbb{R}^2\times(0,T]$, satisfying \eqref{priori_condition} and the assumptions of Theorem \ref{thm1.1}. Then there exists a positive constant $\varepsilon_0=\varepsilon_0(\bar\rho, M)$ depending only on $\mu$, $\lambda$, $\kappa$, $\nu$, $\gamma$, $a$, $\bar\rho$, $\beta$, $N_0$ and $M$ such that
\begin{align}\label{rho_infty_bound}
\sup_{0\leq t\leq T}\|\rho(t)\|_{L^\infty}\leq\frac{7}{4}\bar\rho,
\end{align}
provided $E_0\leq\varepsilon_0$.
\end{lemma}
\begin{proof}
Following an argument analogous to Lemma \(3.7\) in \cite{Li_Xin}, we combine the continuity equation \eqref{model}\(_1\), the effective viscous flux $F$ defined by \eqref{flux_curl}, the estimates \eqref{Gagliardo1}, \eqref{nabla_F_L^p}, \eqref{energy_estimate}, \eqref{A_1_and_A_2}, \eqref{decay_estimate1}, \eqref{u^I_u^P}, \eqref{nabla_F_L^p1}, along with Lemmas \ref{weight_Poincare}, \ref{lem_rho_int_low_bound} and \ref{Zlotnik}, to derive \eqref{rho_infty_bound}.
\end{proof}

\begin{proof}[Proof of Proposition \ref{prop3.1}]
It follows from \eqref{A_1_and_A_2} that, if
$$
E_0\leq\varepsilon_3\triangleq(C(\bar\rho))^{-2},
$$
then
\begin{align}\label{A_1_and_A_2_P}
A_1(T)+A_2(T)+\int_0^T\sigma\|P\|_{L^2}^2dt\leq E_0^{\frac{1}{2}}.
\end{align}
Setting~$\varepsilon\triangleq\min\{\varepsilon_0, \varepsilon_3\}$,
we directly obtain \eqref{uniform_bound} by applying \eqref{rho_infty_bound} and \eqref{A_1_and_A_2_P}. This completes the proof of Proposition \ref{prop3.1}.
\end{proof}

\section{Higher Order Estimates}\label{section4}

Now, for smooth initial data $(\rho_0, u_0, \omega_0)$ satisfying \eqref{initial1} and \eqref{initial2}, we assume that $(\rho, u, \omega)$ is a solution to system \eqref{model}$-$\eqref{initial} on $\mathbb{R}^2\times(0,T]$ satisfying \eqref{priori_condition}. We next derive the higher order regularity estimates of the solution.
\begin{lemma}\label{lem_rho_H^1_W^1}
There exists a positive constant $C$ depending only on $T$, $\mu$, $\lambda$, $\kappa$, $\nu$, $\gamma$, $a$, $\bar\rho$, $\beta$, $N_0$, $M$, $q$, $\|\nabla u_0\|_{L^2}$ and $\|\rho_0\|_{H^1\cap W^{1,q}}$ such that
\begin{equation}\label{rho_H^1_W^1}
	\begin{aligned}[b]
&\sup_{0\leq t\leq T}\bigg(\|\rho\|_{H^1\cap W^{1,q}}+\|\nabla{u}\|_{L^2}^2+\|\nabla\omega\|_{L^2}^2+\|\omega\|_{L^2}^2+t\|\nabla^2u\|_{L^2}^2+t\|\nabla^2\omega\|_{L^2}^2 \bigg)
\\
&+\int_0^T\bigg(\|\nabla^2u\|_{L^2}^2+\|\nabla^2\omega\|_{L^2}^2+\Big(\|\nabla^2u\|_{L^q}+\|\nabla^2\omega\|_{L^q}\Big)^{\frac{1+q}{q}}\bigg)dt
\\
&+\int_0^Tt\Big(\|\nabla^2u\|_{L^q}^2+\|\nabla^2\omega\|_{L^q}^2\Big)dt
\leq C.
	\end{aligned}
\end{equation}
\end{lemma}
\begin{proof}
First, it follows from \eqref{new_B_bound}, \eqref{B_estimate}, Gronwall inequality, and \eqref{energy_estimate} that
\begin{align}\label{new_A_1_estimate}
\sup_{0\leq t\leq T}(\|\nabla u\|_{L^2}^2+\|\nabla{\omega}\|_{L^2}^2+\|\omega\|_{L^2}^2)+\int_0^T\int(\rho|\dot{u}|^2+\rho|\dot{\omega}|^2)dxdt\leq C.
\end{align}
Combining this with \eqref{nabla_u_L^p} and \eqref{nabla_omega_L^p}, one has
\begin{align}\label{nabla_u_nabla_omega}
\int_0^T(\|\nabla{u}\|_{L^4}^4+\|\nabla\omega\|_{L^4}^4)dt\leq C.
\end{align}
Multiplying \eqref{A_2_inequality} by $t$, integrating over $(0,T)$, and then using \eqref{new_A_1_estimate} and \eqref{nabla_u_nabla_omega}, we find
\begin{align}\label{new_A_2_estimate}
\sup_{0\leq t\leq T}t\int(\rho|\dot{u}|^2+\rho|\dot{\omega}|^2)dx+\int_0^Tt(\|\nabla\dot{u}\|_{L^2}^2+\|\nabla\dot\omega\|_{L^2}^2+\|\dot\omega\|_{L^2}^2)dt\leq C.
\end{align}

As in \cite{Huang_Li_Xin}, for $p\in[2,q]$, $|\nabla\rho|^p$ satisfies
\begin{equation}
	\begin{aligned}[b]
&(|\nabla\rho|^p)_t+\Div(u|\nabla\rho|^p)+p|\nabla\rho|^{p-2}(\nabla\rho)^T\nabla u(\nabla\rho)
\\
&+(p-1)|\nabla\rho|^p\Div{u}+p|\nabla\rho|^{p-2}\rho\nabla\rho\cdot\nabla\Div{u}=0.
	\end{aligned}\notag
\end{equation}
Hence,
\begin{equation}\label{rho_L^p_estimate}
	\begin{aligned}[b]
\frac{d}{dt}\|\nabla\rho\|_{L^p}&\leq C\|\nabla{u}\|_{L^\infty}\|\nabla\rho\|_{L^p}+C\|\nabla\Div{u}\|_{L^p}
\\
&\leq C(1+\|\nabla{u}\|_{L^\infty})\|\nabla\rho\|_{L^p}+C\|\rho\dot{u}\|_{L^p},
	\end{aligned}
\end{equation}
where we have used the estimate following from \eqref{nabla_F_L^p}:
\begin{align*}
\|\nabla\Div u\|_{L^p}\leq C(\|\nabla F\|_{L^p}+\|\nabla P\|_{L^p})\leq
C(\|\rho\dot{u}\|_{L^p}+\|\nabla \rho\|_{L^p}).
\end{align*}

Next, by the Gagliardo-Nirenberg inequality, \eqref{new_A_1_estimate}, and \eqref{nabla_F_L^p}, we derive
\begin{equation}
	\begin{aligned}[b]
&\|\Div{u}\|_{L^\infty}+\|\text{curl}~u\|_{L^\infty}\\
\leq& C(\|F\|_{L^\infty}+\|P\|_{L^\infty}+\|\omega\|_{L^\infty}+\|\bar\omega\|_{L^\infty})
\\
\leq& C(q)\|F\|_{L^2}^{\frac{q-2}{2(q-1)}}\|\nabla{F}\|_{L^q}^{\frac{q}{2(q-1)}}+C(q)\|\bar\omega\|_{L^2}^{\frac{q-2}{2(q-1)}}\|\nabla{\bar\omega}\|_{L^q}^{\frac{q}{2(q-1)}}+{C(q)\|\omega\|_{W^{2,q}}}+C(q)
\\
\leq& C(q)\bigg(1+\|\rho\dot{u}\|_{L^q}^{\frac{q}{2(q-1)}}+\|\rho\dot{\omega}\|_{L^q}+{\|\rho\dot{u}\|_{L^2}^{\frac{q-2}{q}}}\bigg).
	\end{aligned}
\end{equation}
The last step comes from the standard $L^q$ elliptic estimate for \eqref{Helmholtz}
\begin{equation}\label{nabla^2_omega}
	\begin{aligned}[b]
\|\omega\|_{W^{2,q}}
&\le C(q)\left(\|\rho\dot\omega\|_{L^q}+\|\bar\omega\|_{L^q}\right)\\
&\le C(q)\left(
\|\rho\dot\omega\|_{L^q}
+\|\bar\omega\|_{L^2}^{\frac{2}{q}}
 \|\nabla\bar\omega\|_{L^2}^{\frac{q-2}{q}}
\right)\\
&\le C(q)\left(
\|\rho\dot\omega\|_{L^q}
+\|\rho\dot u\|_{L^2}^{\frac{q-2}{q}}
\right).
	\end{aligned}
\end{equation}
Moreover, applying the standard $L^q$ elliptic estimate to \eqref{model}$_2$, we get
\begin{align}\label{nabla^2_u}
	\|\nabla^2 u\|_{L^q}
\le C(q)\Big(\|\rho\dot{u}\|_{L^q}+\|\nabla\rho\|_{L^q}+\|\rho\dot\omega\|_{L^q}+\|\rho\dot u\|_{L^2}^{\frac{q-2}{q}}\Big),
\end{align}
which together with Lemma \ref{lem_BKM} yields
\begin{equation}\label{u_infty}
	\begin{aligned}[b]
&\|\nabla{u}\|_{L^\infty}\leq C\bigg(\|\Div{u}\|_{L^\infty}+\|\text{curl}~u\|_{L^\infty}\bigg)\log(e+\|\nabla^2u\|_{L^q})+C\|\nabla{u}\|_{L^2}+C
\\
\leq& C\bigg(1+\|\rho\dot{u}\|_{L^q}+\|\rho\dot{\omega}\|_{L^q}+\|\rho\dot{u}\|_{L^2}^{\frac{q-2}{q}}
\bigg)\log(e+\|\rho\dot{u}\|_{L^q}+\|\rho\dot{\omega}\|_{L^q}+\|\rho\dot{u}\|_{L^2}^{\frac{q-2}{q}}+\|\nabla\rho\|_{L^q}).
	\end{aligned}
\end{equation}
Hence, combining \eqref{rho_L^p_estimate}$-$\eqref{u_infty} gives
\begin{equation}\label{refomulate_rho_L^p_estimate}
\begin{aligned}[b]
&\bigg(\log(e+\|\nabla\rho\|_{L^q})\bigg)_t\leq C(1+\|\nabla u\|_{L^\infty})+C(\|\rho\dot{u}\|_{L^q}+\|\rho\dot\omega\|_{L^q})
\\
\leq& C\bigg(\|\rho\dot{u}\|_{L^q}+\|\rho\dot\omega\|_{L^q}\bigg)+ C\bigg(1+\|\rho\dot{u}\|_{L^q}+\|\rho\dot\omega\|_{L^q}+\|\rho\dot{u}\|_{L^2}^{\frac{q-2}{q}} \bigg)\log(e+\|\nabla\rho\|_{L^q})
\\
&+C\bigg(1+\|\rho\dot{u}\|_{L^q}+\|\rho\dot\omega\|_{L^q}+\|\rho\dot{u}\|_{L^2}^{\frac{q-2}{q}}\bigg)\log(e+\|\rho\dot{u}\|_{L^q}+\|\rho\dot\omega\|_{L^q}+\|\rho\dot{u}\|_{L^2}^{\frac{q-2}{q}}).
\end{aligned}
\end{equation}

In addition, taking $v=\dot{u}$ and $v=\dot{\omega}$, respectively, and $p=q^2$ in \eqref{nabla_F_L^p1}, and then applying H\"{o}lder inequality, we obtain
\begin{equation}\label{rho_dotu_dotomega}
	\begin{aligned}[b]
\|\rho\dot{u}\|_{L^q}+\|\rho\dot{\omega}\|_{L^q}
\leq&\|\rho\dot{u}\|_{L^2}^{\frac{2(q-1)}{q^2-2}}\|\rho\dot{u}\|_{L^{q^2}}^{\frac{q(q-2)}{q^2-2}}+\|\rho\dot{\omega}\|_{L^2}^{\frac{2(q-1)}{q^2-2}}\|\rho\dot{\omega}\|_{L^{q^2}}^{\frac{q(q-2)}{q^2-2}}
\\
\leq &C\|\rho\dot{u}\|_{L^2}^{\frac{2(q-1)}{q^2-2}}\bigg(\|\rho^{\frac{1}{2}}\dot{u}\|_{L^2}+\|\nabla\dot{u}\|_{L^2}\bigg)^{\frac{q(q-2)}{q^2-2}}
\\
&+C\|\rho\dot{\omega}\|_{L^2}^{\frac{2(q-1)}{q^2-2}}\bigg(\|\rho^{\frac{1}{2}}\dot{\omega}\|_{L^2}+\|\nabla\dot{\omega}\|_{L^2}\bigg)^{\frac{q(q-2)}{q^2-2}}
\\
\leq& C\bigg(\|\rho^\frac{1}{2}\dot{u}\|_{L^2}+\|\rho^\frac{1}{2}\dot{\omega}\|_{L^2}\bigg)^{\frac{2(q-1)}{q^2-2}}\bigg(\|\nabla\dot{u}\|_{L^2}+\|\nabla\dot{\omega}\|_{L^2}\bigg)^{\frac{q(q-2)}{q^2-2}}
\\
&+C\bigg(\|\rho^\frac{1}{2}\dot{u}\|_{L^2}+\|\rho^\frac{1}{2}\dot{\omega}\|_{L^2}\bigg).
	\end{aligned}
\end{equation}
Therefore, it follows from \eqref{new_A_1_estimate} and \eqref{new_A_2_estimate} that
\begin{equation}\label{t1+1/q}
	\begin{aligned}[b]
&\int_0^Tt\bigg(\|\rho\dot{u}\|_{L^q}^2+\|\rho\dot{\omega}\|_{L^q}^2\bigg)dt\\
\leq& C\int_0^Tt\bigg(\|\rho^\frac{1}{2}\dot{u}\|_{L^2}^2+\|\rho^\frac{1}{2}\dot{\omega}\|_{L^2}^2\bigg)dt
\\
&+C\int_0^T\bigg(t\Big(\|\rho^\frac{1}{2}\dot{u}\|_{L^2}^2+\|\rho^\frac{1}{2}\dot{\omega}\|_{L^2}^2\Big)\bigg)^{\frac{2(q-1)}{q^2-2}}\bigg(t\Big(\|\nabla\dot{u}\|_{L^2}^2+\|\nabla\dot{\omega}\|_{L^2}^2\Big)\bigg)^{\frac{q(q-2)}{q^2-2}}dt
\\
\leq& C.
	\end{aligned}
\end{equation}
Furthermore, using \eqref{rho_dotu_dotomega} and $0<\frac{q(q+1)(q-2)}{(q^2-2)(q-1)}<1$, we attain
\begin{equation}\label{1+1/q}
	\begin{aligned}[b]
&\int_0^T\bigg(\|\rho\dot{u}\|_{L^q}+\|\rho\dot{\omega}\|_{L^q}\bigg)^{\frac{1+q}{q}}dt\\
\leq& C\int_0^T\bigg(\|\rho^\frac{1}{2}\dot{u}\|_{L^2}^2+\|\rho^\frac{1}{2}\dot{\omega}\|_{L^2}^2\bigg)^{\frac{q+1}{2q}}dt
\\
&+C\int_0^T\bigg(\|\rho^\frac{1}{2}\dot{u}\|_{L^2}^2+\|\rho^\frac{1}{2}\dot{\omega}\|_{L^2}^2\bigg)^{\frac{q^2-1}{q(q^2-2)}}\bigg(\|\nabla\dot{u}\|_{L^2}^2+\|\nabla\dot{\omega}\|_{L^2}^2\bigg)^{\frac{q^2-q-2}{2(q^2-2)}}dt
\\
\leq& C+C\bigg(\int_0^T\Big(\|\rho^\frac{1}{2}\dot{u}\|_{L^2}^2+\|\rho^\frac{1}{2}\dot{\omega}\|_{L^2}^2\Big)dt\bigg)^{\frac{q^2-1}{q(q^2-2)}}
\\
&\times\bigg(\int_0^Tt\Big(\|\nabla\dot{u}\|_{L^2}^2+\|\nabla\dot\omega\|_{L^2}^2\Big)dt\bigg)^{\frac{(q+1)(q-2)}{2(q^2-2)}}\bigg(\int_0^Tt^{-\frac{q(q+1)(q-2)}{(q^2-2)(q-1)}}dt\bigg)^{\frac{q-1}{2q}}
\\
\leq& C.
	\end{aligned}
\end{equation}
Combining \eqref{t1+1/q} and \eqref{1+1/q}, one further derives
\begin{align}\label{t_1+1/q}
\int_0^Tt\bigg(\|\rho\dot{u}\|_{L^q}^2+\|\rho\dot{\omega}\|_{L^q}^2\bigg)dt+\int_0^T\bigg(\|\rho\dot{u}\|_{L^q}+\|\rho\dot{\omega}\|_{L^q}\bigg)^{\frac{1+q}{q}}dt\leq C.
\end{align}
Applying Gronwall inequality to \eqref{refomulate_rho_L^p_estimate} and using \eqref{t_1+1/q}, we deduce that
\begin{align}\label{sup_nabal_rho_L^p}
\sup_{0\leq t\leq T}\|\nabla\rho\|_{L^q}\leq C,
\end{align}
which together with \eqref{u_infty} gives
\begin{align}\label{int_nabla_u_infty}
\int_0^T\|\nabla{u}\|_{L^\infty}dt\leq C.
\end{align}
Taking $p=2$ in \eqref{rho_L^p_estimate} and using Gronwall inequality, \eqref{new_A_1_estimate} and \eqref{int_nabla_u_infty}, we have
\begin{align}\label{sup_nabla_rho_L^2}
\sup_{0\leq t\leq T}\|\nabla\rho\|_{L^2}\leq C.
\end{align}

On the other hand, the standard elliptic estimates, together with \eqref{new_A_1_estimate}, \eqref{new_A_2_estimate}, \eqref{sup_nabla_rho_L^2}, and \eqref{t_1+1/q} give
\begin{align}
\int_0^T\bigg(\|\nabla^2u\|_{L^2}^2+\|\nabla^2\omega\|_{L^2}^2\bigg)dt\leq C,\notag
\\
\sup_{0\leq t\leq T}t\bigg(\|\nabla^2u\|_{L^2}^2+\|\nabla^2\omega\|_{L^2}^2\bigg)\leq C,\notag
\end{align}
and
\begin{align}
\int_0^Tt\bigg(\|\nabla^2u\|_{L^q}^2+\|\nabla^2\omega\|_{L^q}^2\bigg)dt\leq C.\notag
\end{align}
Finally, it follows from \eqref{nabla^2_omega}, \eqref{nabla^2_u}, \eqref{t_1+1/q}, and \eqref{sup_nabal_rho_L^p} that
\begin{align}
\int_0^T\bigg(\|\nabla^2{u}\|_{L^q}+\|\nabla^2{\omega}\|_{L^q}\bigg)^{\frac{1+q}{q}}dt\leq C.\notag
\end{align}
This completes the proof of Lemma \ref{lem_rho_H^1_W^1}.
\end{proof}

\begin{lemma}
There exists a positive constant $C$ depending only on $T$, $\mu$, $\lambda$, $\kappa$, $\nu$, $\gamma$, $a$, $\bar\rho$, $\beta$, $N_0$, $M$, $q$, $\|\nabla u_0\|_{L^2}$, and $\|\nabla(\bar{x}^a\rho_0)\|_{L^2\cap L^q}$ such that
	\begin{align}\label{weight_rho_high_estimate}
\sup_{0\leq t\leq T}\|\bar{x}^a\rho\|_{L^1\cap H^1\cap W^{1,q}}\leq C.
	\end{align}
\end{lemma}
\begin{proof}
	The proof is similar to Lemma $4.2$ in \cite{Li_Xin}.
\end{proof}

\begin{lemma}\label{lem_u_t+omega_t_estimate}
There exists a positive constant $C$ depending only on $T$, $\mu$, $\lambda$, $\kappa$, $\nu$,~$\gamma$, $a$, $\bar\rho$, $\beta$, $N_0$, $M$, $q$, $\|\nabla u_0\|_{L^2}$ and $\|\nabla(\bar{x}^a\rho_0)\|_{L^2\cap L^q}$ such that
\begin{equation}\label{u_t+omega_t_estimate}
		\begin{aligned}[b]
&\sup_{0\leq t\leq T}t\Big(\|\rho^{\frac{1}{2}}u_t\|_{L^2}^2+\|\rho^{\frac{1}{2}}\omega_t\|_{L^2}^2+\|\nabla u\|_{H^1}^2+\|\nabla\omega\|_{H^1}^2\Big)
\\
&+\int_0^Tt\bigg(\|\nabla u_t\|_{L^2}^2+\|\nabla \omega_t\|_{L^2}^2+\|\omega_t\|_{L^2}^2\bigg)dt\leq C.
		\end{aligned}
 \end{equation}
\end{lemma}
\begin{proof}
Differentiating \eqref{model}$_2$ and \eqref{model}$_3$ with respect to $t$, respectively, we obtain
\begin{equation}\label{u_tt}
	\begin{aligned}[b]
&\rho u_{tt}+\rho u\cdot\nabla{u_t}-(\mu+\kappa)\Delta{u_t}-(\mu+\lambda-\kappa)\nabla\Div{u_t}
\\
=&-\rho_t(u_t+u\cdot\nabla u)-\rho u_t\cdot\nabla u-\nabla P_t+2\kappa\nabla^\perp\omega_t,
	\end{aligned}
\end{equation}
and
\begin{align}\label{omega_tt}
\rho\omega_{tt}+\rho u\cdot\nabla\omega_t-\nu\Delta\omega_t+4\kappa \omega_t=-\rho_t(\omega_t+u\cdot\nabla\omega)-\rho u_t\cdot\nabla\omega +2\kappa\text{curl}~u_t.
\end{align}
Multiplying \eqref{u_tt} and \eqref{omega_tt} by $u_t$ and $\omega_t$, respectively, and integrating over $\mathbb{R}^2$, we get
\begin{align*}
	\begin{split}
 &\frac{1}{2}\frac{d}{dt}\int(\rho\vert{}u_t\vert{}^2+\rho\vert\omega_t\vert{}^2)dx+(\mu+\kappa)\int\vert{}\nabla u_t\vert{}^2dx
 \\
 &+(\mu+\lambda-\kappa)\int\vert{}\Div u_t\vert{}^2dx
 +\nu\int\vert{}\nabla \omega_t\vert{}^2dx
 +4\kappa\int\vert{}\omega_t\vert{}^2dx
-4\kappa\int\omega_t~\text{curl}~u_t dx
\\
=&-2\int\rho u\cdot\nabla u_t\cdot u_tdx-\int\rho u\cdot\nabla(u\cdot\nabla u\cdot u_t)dx
-\int\rho u_t\cdot\nabla u\cdot u_t dx
\\
&+\int P_t\Div u_tdx-2\int\rho\omega_t u\cdot\nabla\omega_t dx-\int\rho u\cdot\nabla(\omega_t u\cdot\nabla\omega)dx-\int\rho\omega_t u_t\cdot\nabla\omega dx.
	\end{split}
\end{align*}
Through direct calculations, we further derive
\begin{equation}\label{u_tt+omega_tt}
	\begin{aligned}[b]
&\frac{d}{dt}\int(\rho\vert{}u_t\vert{}^2+\rho\vert\omega_t\vert{}^2)dx+\|\nabla u_t\|_{L^2}^2+\|\nabla\omega_t\|_{L^2}^2+\|\omega_t\|_{L^2}^2
\\
\leq& C\int\rho\vert{}u\vert{}\vert{}u_t\vert{}\bigg(\vert{}\nabla u_t\vert{}+\vert{}\nabla u\vert{}^2+\vert{}u\vert{}\vert{}\nabla^2 u\vert{}\bigg)dx+C\int\rho\vert{}u\vert{}^2\vert{}\nabla u\vert{}\vert{}\nabla u_t\vert{}dx
\\
&+C\int\rho\vert{}u_t\vert{}^2\vert{}\nabla u\vert{}dx + C(\delta)\Vert{}P_t\Vert{}_{L^2}^2+\delta\Vert{}\nabla u_t\Vert{}_{L^2}^2+C\int\rho \vert{}u\vert{}\vert{}\omega_t\vert{}\vert{}\nabla \omega_t\vert{}dx
\\
&+C\int\rho\vert{}u\vert{}\bigg(\vert{}u\vert{}\vert{}\nabla\omega_t\vert{} \vert{}\nabla\omega\vert{}+\vert{}\omega_t\vert{}\vert{}\nabla u\vert{}\vert{}\nabla\omega\vert{}+\vert{}\omega_t\vert{}\vert{}u\vert{}\vert{}\nabla^2\omega\vert{}\bigg)dx
\\
&+C\int\rho\vert{}u_t\vert{}\vert{}\omega_t\vert{}\vert{}\nabla\omega\vert{}dx\\
\triangleq&\sum_{i=1}^8{II}_{i}.
\end{aligned}
\end{equation}

Next, we estimate ${II}_i$. First, using \eqref{weight_L}, \eqref{mass_time_bound}, \eqref{rho_int_low_bound}, and \eqref{wp}, we have, for any $\eta\in(0,1]$ and $s>2$,
\begin{align}
\|u\bar{x}^{-\eta}\|_{L^{\frac{s}{\eta}}}\leq  C(\eta,s).\notag
\end{align}
Then, in view of \eqref{weight_rho_high_estimate}, for any $\eta\in(0,1]$ and $s>2$, we obtain
\begin{align}\label{rho_eta_u}
\|\rho^\eta u\|_{L^{\frac{s}{\eta}}}+\|u\bar{x}^{-\eta}\|_{L^\frac{s}{\eta}}\leq C(\eta,s).
\end{align}
Moreover, the combination of \eqref{wp} and \eqref{mass_time_bound} with
\eqref{rho_int_low_bound} gives for $v=u$ or $v=\omega$,
\begin{align}\label{rho_u_t_L^6}
\|\rho^{\frac{1}{2}}v_t\|_{L^6}\leq C\|\rho^{\frac{1}{2}}v_t\|_{L^2}+C\|\nabla v_t\|_{L^2},
\end{align}

For ${II}_i~(i=1,\cdots,5)$, employing the same arguments in \cite{Li_Xin},
we use \eqref{new_A_1_estimate}, \eqref{rho_eta_u}, \eqref{rho_u_t_L^6}, and \eqref{weight_rho_high_estimate} to derive
\begin{align}\label{u_P_part_estimate}
\sum_{i=1}^5{II}_i\leq\delta\|\nabla u_t\|_{L^2}^2+C(\delta)\left(\|\nabla^2 u\|_{L^2}^2 +\|\rho^{\frac{1}{2}} u_t\|_{L^2}^2+1\right)+C.
\end{align}

Analogously, we estimate ${II}_i~(i=6,7,8)$ as follows:
\begin{equation}
\begin{aligned}[b]
{II}_7=&C\int\rho\vert{}u\vert{}\bigg(\vert{}u\vert{}\vert{}\nabla\omega_t\vert{} \vert{}\nabla\omega\vert{}+\vert{}\omega_t\vert{}\vert{}\nabla u\vert{}\vert{}\nabla\omega\vert{}+\vert{}\omega_t\vert{}\vert{}u\vert{}\vert{}\nabla^2\omega\vert{}\bigg)dx\\
\leq& C\|\rho^{\frac{1}{2}}u\|_{L^8}^2\|\nabla\omega\|_{L^4}\|\nabla\omega_t\|_{L^2}
\\
&+C\|\rho^{\frac{1}{2}}u\|_{L^6}\|\rho^{\frac{1}{2}}\omega_t\|_{L^3}\|\nabla u\|_{L^4}\|\nabla\omega\|_{L^4}+C\|\rho^{\frac{1}{4}}u\|_{L^{12}}^2\|\rho^{\frac{1}{2}}\omega_t\|_{L^3}\|\nabla^2\omega\|_{L^2}
\\
\leq&\delta\|\nabla\omega_t\|_{L^2}^2+C(\delta)\Big(\|\nabla^2u\|_{L^2}^2+\|\nabla^2\omega\|_{L^2}^2+\|\rho^{\frac{1}{2}}\omega_t\|_{L^2}^2+1\ \Big),
\end{aligned}
\end{equation}
and
\begin{equation}\label{omega_part_estimate}
	\begin{aligned}[b]
{II}_6+{II}_8=&C\int\rho\vert{}u\vert{}\vert{}\omega_t\vert{}\vert{}\nabla \omega_t\vert{}dx+C\int\rho\vert{}u_t\vert{}\vert{}\omega_t\vert{}\vert{}\nabla\omega\vert{}dx
\\
\leq &C\|\rho^{\frac{1}{2}}u\|_{L^6}\|\rho^{\frac{1}{2}}\omega_t\|_{L^3}\|\nabla\omega_t\|_{L^2}+C\|\rho^{\frac{1}{2}}u_t\|_{L^6}\|\rho^{\frac{1}{2}}\omega_t\|_{L^3}\|\nabla\omega\|_{L^2}
\\
\leq&\delta\Big(\|\nabla u_t\|_{L^2}^2+\|\nabla\omega_t\|_{L^2}^2\Big)+C(\delta)\Big(1+\|\rho^{\frac{1}{2}}\omega_t\|_{L^2}^2+\|\rho^{\frac{1}{2}}u_t\|_{L^2}^2\Big).
	\end{aligned}
\end{equation}

Substituting \eqref{u_P_part_estimate}$-$\eqref{omega_part_estimate} into \eqref{u_tt+omega_tt}, choosing $\delta$ sufficiently small, and utilizing \eqref{rho_H^1_W^1}, we arrive at
\begin{equation}\label{rho_u_t_omega_t}
	\begin{aligned}[b]
&\frac{d}{dt}\int(\rho\vert{}u_t\vert{}^2+\rho\vert\omega_t\vert{}^2)dx+\|\nabla u_t\|_{L^2}^2+\|\nabla\omega_t\|_{L^2}^2+\|\omega_t\|_{L^2}^2
\\
\leq& C\Big(1+\|\rho^{\frac{1}{2}}\omega_t\|_{L^2}^2+\|\rho^{\frac{1}{2}}u_t\|_{L^2}^2+\|\rho^{\frac{1}{2}}\dot{u}\|_{L^2}^2+\|\rho^{\frac{1}{2}}\dot{\omega}\|_{L^2}^2\Big).
	\end{aligned}
\end{equation}
Here, the standard elliptic estimates, together with \eqref{new_A_1_estimate} and \eqref{sup_nabla_rho_L^2} give
\begin{equation}\label{elliptic_estimate_u_omega}
	\begin{cases}
\|\nabla^2u\|_{L^2}\leq C(\|\rho\dot{u}\|_{L^2}+\|\nabla P\|_{L^2}+\|\nabla\omega\|_{L^2})\leq C(1+\|\rho\dot{u}\|_{L^2}),
\\
\|\nabla^2\omega\|_{L^2}\leq C(\|\rho\dot{\omega}\|_{L^2}+\|\nabla u\|_{L^2}+\|\omega\|_{L^2})\leq C(1+\|\rho\dot{\omega}\|_{L^2}).
	\end{cases}\notag
\end{equation}
Furthermore, using (4.35) and (4.29), we derive
\begin{equation}\label{nabla_u_H^1}
	\begin{aligned}[b]
&\|\nabla u\|_{H^1}+\|\nabla\omega\|_{H^1}+\|\rho^{\frac{1}{2}}u_t\|_{L^2}+\|\rho^{\frac{1}{2}}\omega_t\|_{L^2}
\\
\leq& C(\|\rho^{\frac{1}{2}}\dot{u}\|_{L^2}+\|\rho^{\frac{1}{2}}u\cdot\nabla{u}\|_{L^2}+\|\rho^{\frac{1}{2}}\dot{\omega}\|_{L^2}+\|\rho^{\frac{1}{2}}u\cdot\nabla{\omega}\|_{L^2}+1)
\\
\leq& C(\|\rho^{\frac{1}{2}}\dot{u}\|_{L^2}+\|\rho^{\frac{1}{2}}\dot{\omega}\|_{L^2}+\|\rho^{\frac{1}{2}}{u}\|_{L^6}\|\nabla u\|_{L^3}+\|\rho^{\frac{1}{2}}{u}\|_{L^6}\|\nabla\omega\|_{L^3}+1)
\\
\leq& C(\|\rho^{\frac{1}{2}}\dot{u}\|_{L^2}+\|\rho^{\frac{1}{2}}\dot{\omega}\|_{L^2}+1).
	\end{aligned}
\end{equation}
This, combined with \eqref{rho_u_t_omega_t}, \eqref{new_A_1_estimate} and Gronwall  inequality,  gives \eqref{u_t+omega_t_estimate}. This completes the proof of Lemma \ref{lem_u_t+omega_t_estimate}.
\end{proof}

Throughout the following analysis, we assume that $(\rho, u, \omega)$ is a smooth solution to system \eqref{model}$-$\eqref{initial} on $\mathbb{R}^2\times(0,T]$, satisfying \eqref{priori_condition}, and that the smooth initial data $(\rho_0, u_0, \omega_0)$ satisfy \eqref{initial1}, \eqref{initial2}, \eqref{ad_initial}, and \eqref{compatibility}. In addition to $T$, $\mu$, $\lambda$, $\kappa$, $\nu$, $\gamma$, $a$, $\bar\rho$, $\beta$, $N_0$, $M$, $q$, $\|\nabla u_0\|_{L^2}$, and $\|\nabla(\bar{x}^a\rho_0)\|_{L^2\cap L^q}$, the generic positive constant $C$ may also depend on $\|\nabla^2u_0\|_{L^2}$, $\|\nabla^2\omega_0\|_{L^2}$, $\|\bar{x}^{\delta_0}\nabla^2\rho_0\|_{L^2}$, $\|\bar{x}^{\delta_0}\nabla^2 P(\rho_0)\|_{L^2}$, $\|\nabla^2\rho_0\|_{L^2}$, $\|\nabla^2P(\rho_0)\|_{L^2}$, $\|\nabla^2\rho_0\|_{L^q}$, $\|\nabla^2P(\rho_0)\|_{L^q}$, $\delta_0$, and $\|g_i\|_{L^2}$, where $g_i$ are defined by the compatibility conditions \eqref{compatibility}.

\begin{lemma}\label{lem_rho_u_t_omega_t_H^1}
	The following estimates hold
	\begin{equation}\label{rho_u_t_omega_t_H^1}
		\begin{aligned}[b]
&\sup_{0\leq t\leq T}\Big(\|\rho^{\frac{1}{2}}u_t\|_{L^2}^2+\|\rho^{\frac{1}{2}}\omega_t\|_{L^2}^2+\|\nabla u\|_{H^1}^2+\|\nabla\omega\|_{H^1}^2\Big)
\\
&+\int_0^T\bigg(\|\nabla u_t\|_{L^2}^2+\|\nabla \omega_t\|_{L^2}^2+\|\omega_t\|_{L^2}^2\bigg)dt\leq C.
		\end{aligned}
	\end{equation}
\end{lemma}
\begin{proof}

By virtue of the compatibility conditions \eqref{compatibility}, we can introduce
\begin{align}\label{g_i}
\rho_0^{\frac{1}{2}}\dot{u}_0(x)=g_1,\quad \rho_0^{\frac{1}{2}}\dot{\omega}_0(x)=g_2.
\end{align}
Integrating \eqref{A_2_inequality} over $(0,T)$, together with \eqref{g_i} and \eqref{nabla_u_nabla_omega}, we derive directly that
\begin{align}\label{4.39}
\sup_{0\leq t\leq T}\Big(\|\rho^{\frac{1}{2}}\dot{u}\|_{L^2}^2+\|\rho^{\frac{1}{2}}\dot{\omega}\|_{L^2}^2\Big)+\int_0^T\Big(\|\nabla\dot{u}\|_{L^2}^2+\|\nabla\dot{\omega}\|_{L^2}^2+\|\dot{\omega}\|_{L^2}^2\Big)dt\leq C,
\end{align}
which, together with \eqref{nabla_u_H^1} and \eqref{rho_u_t_omega_t}, gives \eqref{rho_u_t_omega_t_H^1}. This completes the proof of Lemma \ref{lem_rho_u_t_omega_t_H^1}.
\end{proof}
To ensure that the local classical solution can be extended to a global classical solution, we also need to establish the following higher order estimates. The proofs of these estimates are similar to those in \cite{Li_Liang}.
\begin{lemma}\label{lem_noweight_estimate}
It holds that
\begin{align}
&\sup_{0 \le t \le T} \left( \|\bar{x}^{\delta_0}\nabla^2\rho\|_{L^2} + \|\bar{x}^{\delta_0}\nabla^2 P\|_{L^2} \right) \le C, \label{weight_nabla^2_rho}
\\[2ex]
&\sup_{0 \le t \le T} t\bigg(\|\nabla u_t\|_{L^2}^2+\|\nabla\omega_t\|_{L^2}^2+\|\omega_t\|_{L^2}^2\bigg)\notag
\\
& + \int_0^T t \left( \|\rho^{\frac{1}{2}} u_{tt}\|_{L^2}^2+\|\rho^{\frac{1}{2}} \omega_{tt}\|_{L^2}^2 + \|\nabla^2 u_t\|_{L^2}^2+\|\nabla^2 \omega_t\|_{L^2}^2 \right) dt \le C, \label{t_nabla_u_t}
\\[2ex]
&\sup_{0 \le t \le T} \left( \|\nabla^2\rho\|_{L^q} + \|\nabla^2 P\|_{L^q} \right) \le C, \label{nabla^2_rho}
\\[2ex]
&\sup_{0 \le t \le T} t \bigg( \|\rho^{\frac{1}{2}}u_{tt}\|_{L^2} +
\|\rho^{\frac{1}{2}}\omega_{tt}\|_{L^2} +\|\nabla^3 u\|_{L^2 \cap L^q}
+\|\nabla^3 \omega\|_{L^2 \cap L^q} \notag
\\
& + \|\nabla u_t\|_{H^1}
+\|\omega_t\|_{H^2}
 + \|\nabla^2(\rho u)\|_{L^{\frac{q+2}{2}}}
+\|\nabla^2(\rho \omega)\|_{L^{\frac{q+2}{2}}}
\bigg) \notag
\\
&+ \int_0^T t^2 \Big( \|\nabla u_{tt}\|_{L^2}^2+\|\nabla \omega_{tt}\|_{L^2}^2 +\|\omega_{tt}\|_{L^2}^2+ \|\bar{x}^{-1}u_{tt}\|_{L^2}^2 +\|\bar{x}^{-1}\omega_{tt}\|_{L^2}^2\Big) dt \le C. \label{u_tt_omega_tt}
\end{align}
\end{lemma}
\begin{proof}
Based on the similar argument in \cite{Li_Liang} and the estimates \eqref{rho_u_t_omega_t_H^1}, \eqref{4.39}, and \eqref{weight_rho_high_estimate}, we obtain, for $\eta>0$,
\begin{align}\label{weight_u_omega_rho1}
\|\rho^\eta u\|_{L^\infty}+\|\bar{x}^{-\eta}u\|_{L^\infty}+\|\rho^\eta \omega\|_{L^\infty}+\|\bar{x}^{-\eta}\omega\|_{L^\infty}\leq C(\eta).
\end{align}
For $2\leq r\leq q$,  it follows from \eqref{weight_rho_high_estimate}, \eqref{rho_u_t_omega_t_H^1}, and \eqref{weight_u_omega_rho1} that
\begin{align}\label{weight_rho_t_P_t}
\|(\bar{x}^{\frac{a+1}{2}}+|u|+|\omega|)\rho_t\|_{L^r}+\|(\bar{x}^{\frac{a+1}{2}}+|u|+|\omega|)P_t\|_{L^r}\leq C.
\end{align}
Furthermore, \eqref{weight_L}, \eqref{rho_int_low_bound}, \eqref{rho_u_t_omega_t_H^1}, and \eqref{weight_u_omega_rho1} give, for $\eta>0$, and $p>\max\{2,\frac{2}{\eta}\}$
\begin{equation}\label{weight_u_t_omega_t1}
	\begin{cases}
\|\bar{x}^{-\eta}u_t\|_{L^p}+\|\bar{x}^{-\eta}\dot{u}\|_{L^p}\leq C(\eta,p)(1+\|\nabla u_t\|_{L^2}),
\\
\|\bar{x}^{-\eta}\omega_t\|_{L^p}+\|\bar{x}^{-\eta}\dot{\omega}\|_{L^p}\leq C(\eta,p)(1+\|\nabla\omega_t\|_{L^2}).
	\end{cases}
\end{equation}

Using \eqref{rho_H^1_W^1}, \eqref{weight_rho_high_estimate}, \eqref{rho_u_t_omega_t_H^1} and \eqref{weight_u_omega_rho1}$-$\eqref{weight_u_t_omega_t1}, we derive \eqref{weight_nabla^2_rho} by the same argument as in Lemma $4.1$ of \cite{Li_Liang}.  Differentiating  \eqref{model}$_2$ and \eqref{model}$_3$ in space and applying standard elliptic estimates, \eqref{weight_rho_high_estimate}, \eqref{rho_u_t_omega_t_H^1}, \eqref{weight_nabla^2_rho}, \eqref{weight_u_omega_rho1}, and \eqref{weight_u_t_omega_t1}, we obtain
\begin{equation}\label{nabla^3_u_omega_L^2}
	\begin{aligned}[b]
\|\nabla^3u\|_{L^2}+\|\nabla^3\omega\|_{L^2}&\leq C\left(\|\nabla(\rho\dot{u})\|_{L^2}+\|\nabla(\rho\dot{\omega})\|_{L^2}+\|\nabla^2P\|_{L^2}+1\right)
\\
&\leq C(\|\nabla u_t\|_{L^2}+\|\nabla \omega_t\|_{L^2}+\|\omega_t\|_{L^2}+1).
	\end{aligned}
\end{equation}
In addition,  \eqref{weight_rho_high_estimate}, \eqref{rho_u_t_omega_t_H^1}, \eqref{weight_nabla^2_rho}, and \eqref{weight_u_omega_rho1}$-$\eqref{weight_u_t_omega_t1} give
\begin{align}\label{P_t_tt}
\|P_t\|_{H^1}\leq C, \quad \|P_{tt}\|_{L^2}\leq C(1+\|\nabla u_t\|_{L^2}).
\end{align}
Next, multiplying \eqref{u_tt} and \eqref{omega_tt} by $u_{tt}$ and $\omega_{tt}$, respectively,  integrating over space, and
adding the resulting identities, we derive
\begin{equation}\label{estimate_t_nabla_u_t1}
\begin{aligned}[b]
&\frac{d}{dt}B_1(t)
+\int\rho|u_{tt}|^2dx+\int\rho|\omega_{tt}|^2dx\\
=&\int(\rho u)_t\cdot\nabla u_t\cdot\dot{u}dx+\int\rho u\cdot\nabla u_t\cdot\dot{u}_tdx
-2\int\rho u\cdot\nabla u_t\cdot u_{tt}dx\\
&-\int\rho u\cdot\nabla(u\cdot\nabla u)\cdot u_{tt}dx
-\int\rho u_t\cdot\nabla u\cdot u_{tt}dx
-\int P_{tt}\Div u_{t}dx\\
&-2\int\rho u\cdot\nabla\omega_t\omega_{tt}dx
 -\int\rho u\cdot\nabla(u\cdot\nabla\omega)\omega_{tt}dx
-\int\rho \omega_{tt}u_t\cdot\nabla\omega dx\\
&+\int(\rho u)_t\cdot\nabla\omega_t\cdot\dot{\omega}dx+\int\rho u\cdot\nabla\omega_t\dot{\omega}_tdx
\\
\triangleq&\sum_{i=1}^{11}{III}_i,
\end{aligned}
\end{equation}
where
\begin{equation}\label{4.45}
	\begin{aligned}[b]
B_1(t)&=\frac{\mu+\kappa}{2}\|\nabla u_t\|_{L^2}^2+\frac{\mu+\lambda-\kappa}{2}\|\Div u_t\|_{L^2}^2+\frac{\nu}{2}\|\nabla\omega_t\|_{L^2}^2+2\kappa\|\omega_t\|_{L^2}^2
\\
&\quad+\int\rho u\cdot\nabla u_t\cdot\dot udx
+\int\rho u\cdot\nabla\omega_t\dot\omega dx-\int P_t\Div u_tdx
-2\kappa\int\omega_t\text{curl}~u_tdx.
	\end{aligned}
\end{equation}
As in \cite{Li_Liang}, we estimate ${III}_i~(i=1,...,6)$ by using \eqref{rho_u_t_omega_t_H^1}, \eqref{weight_u_omega_rho1}$-$\eqref{weight_u_t_omega_t1} and obtain
\begin{align}\label{NS_parts}
	\begin{split}
\sum_{i=1}^6\left|{III}_i\right|&\leq\delta\|\rho^{\frac{1}{2}}u_{tt}\|_{L^2}^2+\varepsilon\|\nabla^2u_t\|_{L^2}^2+C(\delta,\varepsilon)\big(1+\|\nabla u_t\|_{L^2}^2+\|\nabla u_t\|_{L^2}^4\big).
	\end{split}
\end{align}
For the microrotation terms,   \eqref{wp}, \eqref{rho_int_low_bound}, \eqref{weight_rho_high_estimate}, and \eqref{rho_u_t_omega_t_H^1} give
\begin{align*}
\|\rho^{\frac{1}{4}}u_t\|_{L^4}+\|\rho^{\frac{1}{4}}\omega_t\|_{L^4}\leq C(1+\|\nabla u_t\|_{L^2}+\|\nabla\omega_t\|_{L^2}).
\end{align*}
Hence,
\begin{equation}\label{new_parts1}
	\begin{aligned}[b]
\sum_{i=7}^9\left|{III}_i\right|\leq& C\big(\|\rho^{\frac{1}{2}}u\|_{L^\infty}\|\rho^{\frac{1}{2}}\omega_{tt}\|_{L^2}\|\nabla\omega_t\|_{L^2}+\|\rho^{\frac{1}{2}}u_t\|_{L^4}\|\rho^{\frac{1}{2}}\omega_{tt}\|_{L^2}\|\nabla\omega\|_{L^4}\big)
\\
& +\delta\|\rho^{\frac{1}{2}}\omega_{tt}\|_{L^2}^2+C(\delta)\int\rho|u|^2(|\nabla u|^2|\nabla\omega|^2+|u|^2|\nabla^2\omega|^2)dx
\\
\leq&\delta\|\rho^{\frac{1}{2}}\omega_{tt}\|_{L^2}^2+C(\delta)(\|\nabla u_{t}\|_{L^2}^2+\|\nabla\omega_{t}\|_{L^2}^2+\|\omega_t\|_{L^2}^2+1).
	\end{aligned}
\end{equation}
For $III_{10}$ and $III_{11}$, we note that  \eqref{weight_rho_high_estimate}, and \eqref{weight_u_omega_rho1}$-$\eqref{weight_u_t_omega_t1} imply
\begin{align}\label{bar_rho_u_t}
\|\bar{x}(\rho u)_t\|_{L^q}\leq C(\|\nabla u_t\|_{L^2}+1),
\end{align}
which, together with  \eqref{Gagliardo1}, \eqref{weight_rho_high_estimate}, \eqref{rho_u_t_omega_t_H^1} and \eqref{weight_u_omega_rho1}$-$\eqref{weight_u_t_omega_t1}, yields
\begin{equation}\label{new_parts2}
	\begin{aligned}[b]
&\big|{III}_{10}\big|+\big|{III}_{11}\big|
\\
\leq& \|\rho^{\frac{1}{2}}u\|_{L^\infty}\|\nabla\omega_t\|_{L^2}\|\rho^{\frac{1}{2}}\omega_{tt}\|_{L^2}+\|\rho^\frac{1}{2}u\|_{L^\infty}\|\nabla\omega_t\|_{L^4}\|\rho^{\frac{1}{2}}u_t\|_{L^4}\|\nabla\omega\|_{L^2}
\\
&+\|\bar x(\rho u)_t\|_{L^q}\|\nabla\omega_t\|_{L^2}\|\bar x^{-1}\dot\omega\|_{L^{\frac{2q}{q-2}}}
+\|\rho^{\frac12}u\|_{L^\infty}^2\|\nabla\omega_t\|_{L^2}^2
\\
\leq&\delta\|\rho^{\frac{1}{2}}\omega_{tt}\|_{L^2}^2+C(\delta)\|\nabla\omega_t\|_{L^2}^2+C(1+\|\nabla u_t\|_{L^2})\|\nabla\omega_t\|_{L^2}^\frac{1}{2}\|\nabla^2\omega_t\|_{L^2}^\frac{1}{2}
\\
&+C\bigl(1+\|\nabla u_t\|_{L^2}\bigr)\|\nabla\omega_t\|_{L^2}\bigl(1+\|\nabla\omega_t\|_{L^2}\bigr)
\\
\leq& \delta\|\rho^{\frac{1}{2}}\omega_{tt}\|_{L^2}^2+\varepsilon\|\nabla^2\omega_t\|_{L^2}^2+C(\delta,\varepsilon)+C(\delta,\varepsilon)\big(\|\nabla u_t\|_{L^2}^2+\|\nabla\omega_t\|_{L^2}^2+\|\omega_t\|_{L^2}^2\big)^2.
	\end{aligned}
\end{equation}
Putting \eqref{NS_parts}, \eqref{new_parts1}, and \eqref{new_parts2} into \eqref{estimate_t_nabla_u_t1} and choosing $\delta$  sufficiently small, we have
\begin{equation}\label{estimate_t_nabla_u_t2}
	\begin{aligned}[b]
&\frac{d}{dt}B_1(t)
+\frac{1}{2}\int\rho|u_{tt}|^2dx+\frac{1}{2}\int\rho|\omega_{tt}|^2dx
\\
\leq&\varepsilon(\|\nabla^2u_t\|_{L^2}^2+\|\nabla^2\omega_t\|_{L^2}^2)+C(\varepsilon)\Big(1+\big(\|\nabla u_t\|_{L^2}^2+\|\nabla\omega_t\|_{L^2}^2+\|\omega_t\|_{L^2}^2\big)^2\Big).
	\end{aligned}
\end{equation}
In addition, it follows from H\"older inequality and \eqref{4.39} that
\begin{align*}
	\begin{array}{l}
\Big|\displaystyle\int\rho u\cdot\nabla u_t\cdot\dot udx\Big|\leq\frac{\mu}{16}\|\nabla u_t\|_{L^2}^2+C,
	\\
\left|2\kappa\displaystyle\int\omega_t\text{curl}~u_tdx\right|\leq\frac{4\kappa^2}{\mu+2\kappa}\|\omega_t\|_{L^2}^2+\frac{\mu+2\kappa}{4}\|\nabla u_t\|_{L^2}^2,
	\end{array}
	\begin{array}{l}
\Big|\displaystyle\int\rho u\cdot\nabla\omega_t\dot\omega dx\Big|\leq\frac{\nu}{4}\|\nabla\omega_t\|_{L^2}^2+C,
	\\
\Big|\displaystyle\int P_t\Div u_tdx\Big|\leq\frac{\mu}{16}\|\nabla u_t\|_{L^2}^2+C,
	\end{array}
\end{align*}
which, together with \eqref{4.45}, gives,
\begin{equation}\label{bound_B_1_t}
	\begin{aligned}[b]
&\frac{\mu}{8}\|\nabla u_t\|_{L^2}^2+\frac{\nu}{4}\|\nabla\omega_t\|_{L^2}^2+\frac{2\kappa\mu}{\mu+2\kappa}\|\omega_t\|_{L^2}^2-C
\\
\leq &B_1(t)\leq C\big(1+\|\nabla u_t\|_{L^2}^2+\|\nabla\omega_t\|_{L^2}^2+\|\omega_t\|_{L^2}^2\big).
\end{aligned}
\end{equation}
Moreover, using \eqref{u_tt} and \eqref{omega_tt} and the standard elliptic estimates, we get
\begin{equation}\label{elliptic_nabla^2_u_t}
	\begin{aligned}[b]
\|\nabla^2u_t\|_{L^2}^2+\|\nabla^2\omega_t\|_{L^2}^2\leq &C\big(\|\rho^{\frac{1}{2}}u_{tt}\|_{L^2}^2+\|\rho^{\frac{1}{2}}\omega_{tt}\|_{L^2}^2
\\
&+\|\nabla u_t\|_{L^2}^2+\|\nabla\omega_t\|_{L^2}^2+\|\omega_t\|_{L^2}^2+1\big).
\end{aligned}
\end{equation}
Multiplying \eqref{estimate_t_nabla_u_t2} by $t$, selecting $\varepsilon$ small enough, and using  \eqref{rho_u_t_omega_t_H^1}, \eqref{bound_B_1_t}, \eqref{elliptic_nabla^2_u_t} and Gronwall inequality, we derive \eqref{t_nabla_u_t}.

Using  \eqref{rho_H^1_W^1}, \eqref{weight_rho_high_estimate}, \eqref{rho_u_t_omega_t_H^1}, \eqref{t_nabla_u_t},  and \eqref{weight_u_omega_rho1}$-$\eqref{nabla^3_u_omega_L^2}, we can derive \eqref{nabla^2_rho} by the same argument as in Lemma $4.3$ of \cite{Li_Liang}. Applying the $L^q$ elliptic estimates for the spatial derivatives of \eqref{model}$_2$ and \eqref{model}$_3$, combined with
\eqref{Gagliardo1}, \eqref{weight_rho_high_estimate}, \eqref{rho_u_t_omega_t_H^1}, \eqref{nabla^2_rho}, \eqref{weight_u_omega_rho1}, \eqref{weight_u_t_omega_t1}, and \eqref{nabla^3_u_omega_L^2}, we obtain
\begin{equation}\label{nabla^3_u_omega_L^q}
	\begin{aligned}
		[b]
\|\nabla^3u\|_{L^q}+\|\nabla^3\omega\|_{L^q}\leq &C(\|\nabla u_t\|_{L^2}+\|\nabla \omega_t\|_{L^2}+\|\omega_t\|_{L^2}+1)
\\
&+C(\|\nabla u_t\|_{L^q}+\|\nabla \omega_t\|_{L^q}).
	\end{aligned}
\end{equation}

Differentiating \eqref{u_tt} and \eqref{omega_tt} with respect to $t$ leads to
\begin{align}
&\rho u_{ttt}+\rho u\cdot\nabla u_{tt}-(\mu+\kappa)\Delta u_{tt}-(\mu+\lambda-\kappa)\nabla\Div u_{tt}-2\kappa\nabla^\perp\omega_{tt} \notag
 \\
=&-2\rho_tu_{tt}-\rho_{tt}u_t-2(\rho u)_t\cdot\nabla u_t-\rho_{tt}u\cdot\nabla u-2\rho_tu_t\cdot\nabla u-\rho u_{tt}\cdot\nabla u-\nabla P_{tt}, \label{u_ttt}
 \\
&\rho\omega_{ttt}+\rho u\cdot\nabla\omega_{tt}-\nu\Delta\omega_{tt}+4\kappa\omega_{tt}-2\kappa\text{curl}~u_{tt}    \notag
 \\
=&-2\rho_t\omega_{tt}-\rho_{tt}\omega_t-2(\rho u)_t\cdot\nabla\omega_t-\rho_{tt}u\cdot\nabla\omega-2\rho_tu_t\cdot\nabla\omega-\rho u_{tt}\cdot\nabla\omega.  \label{omega_ttt}
\end{align}
Multiplying \eqref{u_ttt} and \eqref{omega_ttt} by $u_{tt}$ and $\omega_{tt}$, respectively,  integrating by parts over $\mathbb{R}^2$, and summing the resulting identities, we get
\begin{equation}\label{rho_u_tt_estimate}
	\begin{aligned}[b]
&\frac{1}{2}\frac{d}{dt}\int\big(\rho|u_{tt}|^2+\rho|\omega_{tt}|^2\big)dx+(\mu+\kappa)\|\nabla u_{tt}\|_{L^2}^2
+(\mu+\lambda-\kappa)\|\Div u_{tt}\|_{L^2}^2
\\
& +\nu\|\nabla\omega_{tt}\|_{L^2}^2+4\kappa\|\omega_{tt}\|_{L^2}^2-4\kappa\int\omega_{tt}\text{curl}~u_{tt}dx
\\
=&2\int\Div(\rho u)|u_{tt}|^2dx+\int\Div(\rho u)_tu_t\cdot u_{tt}dx-2\int(\rho u)_t\cdot\nabla u_t\cdot u_{tt}dx
\\
& -\int\rho_{tt}u\cdot\nabla u\cdot u_{tt}dx
-2\int\rho_tu_t\cdot\nabla u\cdot u_{tt}dx-\int\rho u_{tt}\cdot\nabla u\cdot u_{tt}dx+\int P_{tt}\Div u_{tt}dx
\\
& +2\int\Div(\rho u)|\omega_{tt}|^2dx+\int\Div(\rho u)_t\omega_t\omega_{tt}dx-2\int(\rho u)_t\cdot\nabla\omega_t\omega_{tt}dx
\\
& -\int\rho_{tt}\omega_{tt}u\cdot\nabla\omega dx-2\int\rho_t\omega_{tt}u_t\cdot\nabla\omega dx-\int\rho\omega_{tt} u_{tt}\cdot\nabla\omega dx
\\
\triangleq&\sum_{i=1}^{13}{IV}_i.
	\end{aligned}
\end{equation}
Employing \eqref{rho_u_t_omega_t_H^1}, \eqref{weight_u_omega_rho1}$-$\eqref{weight_u_t_omega_t1}, \eqref{P_t_tt}, and \eqref{bar_rho_u_t}, we estimate ${IV}_i,~i=1, \cdots,7$ as in \cite{Li_Liang}. This gives
\begin{align}\label{iv_1-7}
	\sum_{i=1}^7\left|{IV}_i\right|\leq\delta\|\nabla u_{tt}\|_{L^2}^2+C(\delta)\big(\|\rho^{\frac{1}{2}}u_{tt}\|_{L^2}^2+\|\nabla u_t\|_{L^2}^4+\|\nabla u_t\|_{L^2}^2+1\big).
\end{align}
For  $IV_i,~i=8,\cdots,13$,  we first use \eqref{weight_L} and \eqref{rho_int_low_bound} to obtain
\begin{align}\label{weight_omega_tt}
\|\bar{x}^{-1}\omega_{tt}\|_{L^{\frac{2q}{q-2}}}+\|\bar{x}^{-\frac{1}{2}}\omega_{tt}\|_{L^{\frac{4q}{q-2}}}\leq C(\|\rho^{\frac{1}{2}}\omega_{tt}\|_{L^2}+\|\nabla\omega_{tt}\|_{L^2}).
\end{align}
By \eqref{weight_u_omega_rho1} and Young inequality, we have
\begin{equation}\label{iv_8}
	\begin{aligned}[b]
|IV_8|&\leq C\int\rho|u||\nabla\omega_{tt}||\omega_{tt}|dx\leq C\|\rho^{\frac{1}{2}}u\|_{L^\infty}\|\rho^{\frac{1}{2}}\omega_{tt}\|_{L^2}\|\nabla\omega_{tt}\|_{L^2}
\\
&\leq\delta\|\nabla\omega_{tt}\|_{L^2}^2+C(\delta)\|\rho^{\frac{1}{2}}\omega_{tt}\|_{L^2}^2.
	\end{aligned}
\end{equation}
For $IV_{9}$, $IV_{10}$,  integrating by parts and using  \eqref{bar_rho_u_t}, \eqref{weight_u_t_omega_t1}, \eqref{weight_omega_tt}, and H\"older inequality, we obtain
\begin{equation}
	\begin{aligned}[b]
&|IV_9|+|{IV}_{10}|
\\
\leq& C\|\bar{x}(\rho u)_t\|_{L^q}\|\bar{x}^{-1}\omega_{tt}\|_{L^{\frac{2q}{q-2}}}\|\nabla\omega_t\|_{L^2}+C\|\bar{x}(\rho u)_t\|_{L^q}\|\bar{x}^{-1}\omega_{t}\|_{L^{\frac{2q}{q-2}}}\|\nabla\omega_{tt}\|_{L^2}
\\
\leq&\delta\|\nabla\omega_{tt}\|_{L^2}^2+C(\delta)\|\rho^{\frac{1}{2}}\omega_{tt}\|_{L^2}^2+C(\delta)\big(\|\nabla u_t\|_{L^2}^4+\|\nabla\omega_{t}\|_{L^2}^4+1\big),
	\end{aligned}
\end{equation}
For $IV_{11}$, integrating by parts with $\rho_{tt}=-\Div(\rho u)_t$,  and using \eqref{weight_L}, \eqref{rho_int_low_bound},
\eqref{rho_u_t_omega_t_H^1}, \eqref{weight_u_omega_rho1}, \eqref{bar_rho_u_t}, \eqref{weight_omega_tt} and H\"older inequality, for $0<r<\frac{2}{q}$, we have
\begin{equation}
	\begin{aligned}[b]
|IV_{11}|
\le& C\|\bar{x}(\rho u)_t\|_{L^q}
\bigg(
\|\nabla u\|_{L^4}\|\nabla\omega\|_{L^4}
\|\bar{x}^{-1}\omega_{tt}\|_{L^\frac{2q}{q-2}}\\
&
+\|\bar{x}^{-r}u\|_{L^\infty}
\|\nabla^2\omega\|_{L^2}
\|\bar{x}^{-(1-r)}\omega_{tt}\|_{L^\frac{2q}{q-2}}\\
&
+\|\bar{x}^{-1}u\|_{L^\infty}
\|\nabla\omega\|_{L^\frac{2q}{q-2}}
\|\nabla\omega_{tt}\|_{L^2}
\bigg)
\\
\leq&\delta\|\nabla\omega_{tt}\|_{L^2}^2+C(\delta)\|\rho^{\frac{1}{2}}\omega_{tt}\|_{L^2}^2+C(\delta)\big(\|\nabla u_t\|_{L^2}^2+1 \big).
	\end{aligned}
\end{equation}
By virtue of \eqref{wp}, \eqref{rho_u_t_omega_t_H^1}, \eqref{weight_rho_t_P_t}, \eqref{weight_u_t_omega_t1}, and \eqref{weight_omega_tt}, one has
 \begin{equation}\label{IV_12_13}
	\begin{aligned}[b]
&|IV_{12}|+|IV_{13}|
\\
\leq& C\|\bar{x}\rho_t\|_{L^q}\|\bar{x}^{-\frac{1}{2}}u_t\|_{L^\frac{4q}{q-2}}\|\bar{x}^{-\frac{1}{2}}\omega_{tt}\|_{L^\frac{4q}{q-2}}\|\nabla\omega\|_{L^2}+C\|\rho^{\frac{1}{2}}u_{tt}\|_{L^6}\|\nabla\omega\|_{L^2}\|\rho^{\frac{1}{2}}\omega_{tt}\|_{L^3}
\\
\leq&\delta\big(\|\nabla\omega_{tt}\|_{L^2}^2+\|\nabla u_{tt}\|_{L^2}^2 \big) +C(\delta)\big(\|\rho^{\frac{1}{2}}\omega_{tt}\|_{L^2}^2+\|\rho^{\frac{1}{2}}u_{tt}\|_{L^2}^2\big)
\\
&+C(\delta)\big(\|\nabla u_t\|_{L^2}^2+\|\nabla\omega_{t}\|_{L^2}^2+\|\omega_{t}\|_{L^2}^2+1 \big).
	\end{aligned}
 \end{equation}
Moreover, the coupled dissipation satisfies
\begin{align*}
\begin{split}
&(\mu+\kappa)\|\nabla u_{tt}\|_{L^2}^2
 +(\mu+\lambda-\kappa)\|\Div u_{tt}\|_{L^2}^2
 +\nu\|\nabla\omega_{tt}\|_{L^2}^2
 +4\kappa\|\omega_{tt}\|_{L^2}^2
 -4\kappa\int\omega_{tt}\text{curl}~u_{tt}\,dx\\
=&(2\mu+\lambda)\|\Div u_{tt}\|_{L^2}^2
 +\mu\|\text{curl}~u_{tt}\|_{L^2}^2
 +\kappa\|\text{curl}~u_{tt}-2\omega_{tt}\|_{L^2}^2
 +\nu\|\nabla\omega_{tt}\|_{L^2}^2\\
\ge& C\bigl(\|\nabla u_{tt}\|_{L^2}^2
 +\|\nabla\omega_{tt}\|_{L^2}^2+\|\omega_{tt}\|_{L^2}^2\bigr),
	\end{split}
\end{align*}
Combining this estimate with \eqref{rho_u_tt_estimate}$-$\eqref{IV_12_13} and choosing $\delta$ sufficiently small, we obtain
\begin{equation}\label{rho_u_tt_omega_tt_estimate}
	\begin{aligned}[b]
&\frac{d}{dt}\left(\|\rho^{\frac12}u_{tt}\|_{L^2}^2
       +\|\rho^{\frac12}\omega_{tt}\|_{L^2}^2\right)
 +\left(\|\nabla u_{tt}\|_{L^2}^2
       +\|\nabla\omega_{tt}\|_{L^2}^2+\|\omega_{tt}\|_{L^2}^2\right)\\
\le& C\left(\|\rho^{\frac12}u_{tt}\|_{L^2}^2
       +\|\rho^{\frac12}\omega_{tt}\|_{L^2}^2\right)
 +C\left(1+\|\nabla u_t\|_{L^2}^4+\|\nabla\omega_t\|_{L^2}^4+\|\omega_t\|_{L^2}^4\right).
	\end{aligned}
\end{equation}
Multiplying \eqref{rho_u_tt_omega_tt_estimate} by $t^2$ and using Gronwall inequality and \eqref{t_nabla_u_t}, we obtain
\begin{align}\label{rho_u_tt_omega_tt_estimate1}
\sup_{0\leq t\leq T}t^2\left(\|\rho^{\frac12}u_{tt}\|_{L^2}^2
       +\|\rho^{\frac12}\omega_{tt}\|_{L^2}^2\right)+\int_0^Tt^2\left(\|\nabla u_{tt}\|_{L^2}^2
       +\|\nabla\omega_{tt}\|_{L^2}^2+\|\omega_{tt}\|_{L^2}^2\right)dt\leq C.
\end{align}
In addition, it follows from
\eqref{weight_L}, \eqref{rho_int_low_bound}, and
\eqref{rho_u_tt_omega_tt_estimate1} that
\begin{align}
\int_0^Tt^2\left(\|\bar{x}^{-1}u_{tt}\|_{L^2}^2+\|\bar{x}^{-1}\omega_{tt}\|_{L^2}^2\right)dt\leq C.\notag
\end{align}
Then, based on \eqref{elliptic_nabla^2_u_t}, \eqref{t_nabla_u_t}, \eqref{rho_u_tt_omega_tt_estimate1}, one gets
\begin{align}\label{omega_H^2}
	\sup_{0\leq t\leq T}t\left(\|\nabla u_t\|_{H^1}+\|\omega_t\|_{H^2}\right)\leq C.
\end{align}
In view of  \eqref{Gagliardo1}, \eqref{nabla^3_u_omega_L^2}, \eqref{nabla^3_u_omega_L^q}, and \eqref{omega_H^2},   one gets
\begin{align}\label{4.72}
\sup_{0\leq t\leq T}t\left(\|\nabla^3u\|_{L^2\cap L^q}+\|\nabla^3\omega\|_{L^2\cap L^q}\right)\leq C.
\end{align}
 Similar to \cite{Li_Liang},  we  utilize \eqref{weight_rho_high_estimate}, \eqref{rho_u_t_omega_t_H^1}, \eqref{weight_nabla^2_rho}, \eqref{nabla^2_rho}, \eqref{weight_u_omega_rho1}, and \eqref{4.72} to obtain
\begin{equation}
	\begin{aligned}[b]
t\|\nabla^2(\rho u)\|_{L^\frac{q+2}{2}}\leq&
Ct\|\bar{x}^{\delta_0}\nabla^2\rho\|_{L^2}^{\frac{2}{q+2}}\|\nabla^2\rho\|_{L^q}^\frac{q}{q+2}\|\bar{x}^{-\frac{2\delta_0}{q+2}}u\|_{L^\infty}
\\
&+Ct\|\nabla\rho\|_{L^q}\|\nabla u\|_{L^\frac{q(q+2)}{q-2}}+Ct\|\nabla^2u\|_{L^{\frac{q+2}{2}}}\leq C,
	\end{aligned}\notag
\end{equation}
and, since \eqref{rho_H^1_W^1} and \eqref{rho_u_t_omega_t_H^1} imply $\|\omega\|_{H^2}\leq C$,
\begin{equation}
	\begin{aligned}[b]
t\|\nabla^2(\rho \omega)\|_{L^\frac{q+2}{2}}
\leq &Ct\|\omega\nabla^2\rho\|_{L^{\frac{q+2}{2}}}+Ct\|\nabla\rho\nabla\omega\|_{L^{\frac{q+2}{2}}}+Ct\|\rho{\nabla^2}\omega\|_{L^{\frac{q+2}{2}}}
\\
\leq& Ct\|\omega\|_{L^\infty}\|\nabla^2\rho\|_{L^{\frac{q+2}{2}}}+Ct\|\nabla\rho\|_{L^q}\|\nabla\omega\|_{L^{\frac{q(q+2)}{q-2}}}
+Ct\|\nabla^2\omega\|_{L^{{\frac{q+2}{2}}}}\\
\leq& C.
 	\end{aligned}\notag
\end{equation}
The proof of Lemma \ref{lem_noweight_estimate} is finished.
\end{proof}

\section{Proofs of Theorems \ref{thm1.1} and \ref{thm1.2}}\label{section5}
With all the {\it a priori} estimates established in Sections \ref{section3} and \ref{section4}, we are ready to complete the proofs of Theorems \ref{thm1.1} and \ref{thm1.2}.
\begin{proof}[\noindent\textbf{\textit{Proof of Theorem \ref{thm1.1}}}]

By Lemma \ref{local}, there exists $T_*>0$ such that the Cauchy problem \eqref{model}$-$\eqref{initial} admits a unique strong solution $(\rho, u, \omega)$ on $\mathbb{R}^2\times(0,T_*]$. We will use the {\it a priori} estimates established in Proposition \ref{prop3.1} and Lemmas \ref{lem_rho_H^1_W^1}$-$\ref{lem_u_t+omega_t_estimate} to extend the local strong solution to a global one.

First, it follows from \eqref{A_1}, \eqref{A_2}, and \eqref{initial1}  that
\begin{align*}
A_1(0)+A_2(0)=0, \quad \rho_0\leq\bar\rho.
\end{align*}
Thus, there exists $T_1\in(0,T_*]$ such that \eqref{priori_condition} holds when $T=T_1$.

Next, let
\begin{align}\label{sup_T}
	T^*=\sup\{T~|~\eqref{priori_condition}~\text{holds}\}.
\end{align}
Clearly, $T^*\geq T_1>0$. Hence, for any $0<\tau<T\leq T^*$ with $T$ finite, it follows from \eqref{u_t+omega_t_estimate} that for any $2\leq q<\infty$,
\begin{align}\label{nabla_u_omega_tau}
\nabla u\in C([\tau,T]; L^2\cap L^q), \quad \nabla\omega\in C([\tau,T];  L^2\cap L^q),
\end{align}
where we have used the embedding
\begin{align*}
L^\infty(\tau,T; H^1)\cap H^1(\tau,T; H^{-1})
\hookrightarrow C([\tau,T];L^q), \quad \forall q\in[2,\infty).
\end{align*}
Moreover, one deduces from \eqref{rho_H^1_W^1}, \eqref{weight_rho_high_estimate}, and Lemma $2.3$ in \cite{Lions} that
\begin{align}\label{C_rho}
\rho\in C([0,T];L^1\cap H^1\cap W^{1,q}).
\end{align}

Finally, we claim that
\begin{align}
T^*=\infty.
\end{align}
Otherwise, suppose that $T^*<\infty$. By Proposition \ref{prop3.1}, \eqref{uniform_bound} holds when $T=T^*$. It follows from \eqref{energy_estimate}, \eqref{weight_rho_high_estimate}, \eqref{nabla_u_omega_tau}, and \eqref{C_rho} that $(\rho(x,T^*), u(x,T^*), \omega(x,T^*))$ satisfies all the requirements of the initial conditions in \eqref{initial1}, except for $u(\cdot,T^*)\in\dot{H}^\beta$. The local strong existence assertion of Lemma \ref{local} does not require the $\dot{H}^\beta$ condition. That condition is used in the estimates near the original initial time, and it is therefore unnecessary to impose it again when restarting the
strong solution at $T^*>0$. Thus, by Lemma \ref{local}, there exists a time $T^{**}>T^*$ such that \eqref{priori_condition} still holds when $T=T^{**}$. This contradicts the definition of $T^*$ in \eqref{sup_T}. Hence, $T^*=\infty$ holds. It then follows from Lemma \ref{local} and Lemmas \ref{lem_rho_H^1_W^1}$-$\ref{lem_u_t+omega_t_estimate} that, for any $0<T<T^*=\infty$, $(\rho, u, \omega)$ is indeed the unique strong solution defined on $\mathbb{R}^2\times(0,T]$. This completes the proof of Theorem 1.1.
\end{proof}

\begin{proof}[\noindent\textbf{\textit{Proof of Theorem \ref{thm1.2}}}]
By Lemma \ref{local}, under \eqref{ad_initial} and \eqref{compatibility}, problem \eqref{model}$-$\eqref{initial} admits a unique local classical solution $(\rho,u,\omega)$. By uniqueness, it coincides with the global strong solution obtained in Theorem \ref{thm1.1}. It follows from Lemmas \ref{lem_rho_u_t_omega_t_H^1} and \ref{lem_noweight_estimate}, together with \eqref{weight_L} and \eqref{rho_int_low_bound}, that the norm bounds in \eqref{regularity1} hold on every finite interval of existence.

Moreover, by interpolation and \eqref{rho_u_t_omega_t_H^1}, \eqref{t_nabla_u_t}, \eqref{nabla^3_u_omega_L^2}, and \eqref{nabla^3_u_omega_L^q}, we have $\nabla^3u\in L^1(0,T;L^2\cap L^q)$. This, combined with \eqref{weight_nabla^2_rho}, \eqref{nabla^2_rho}, and the standard continuity argument applied to the equations satisfied by $\nabla^2\rho$ and $\nabla^2P$, gives
\begin{align*}
\nabla^2\rho,\ \nabla^2P\in C([0,T];L^2\cap L^q).
\end{align*}
Therefore, using \eqref{4.39}, Lemma \ref{local}, and the same continuation argument as in the proof of Theorem \ref{thm1.1}, we extend the local classical solution globally, with \eqref{regularity1} holding for every $0<T<\infty$. This completes the proof of Theorem \ref{thm1.2}.
\end{proof}

\vspace{6mm}
\noindent {\bf Acknowledgements:}
The first and second authors were supported by the Basic and Applied Basic Research Foundation of Guangdong Province (No. 2026A1515012081) and the National Natural Science Foundation of China (No. 12571233). The third author was supported by the National Natural Science Foundation of China (No. 12571233) and the Basic and Applied Basic Research Foundation of Guangdong Province (No. 2026A1515012081).

\vspace{6mm}

\noindent {\bf Data Availability:} No datasets were generated or analysed for this study.

\noindent {\bf Conflicts of Interest:} The authors declare that they have no conflicts of interest.

\end{document}